\documentclass[draft]{amsart}

\newcommand{\cl}{\operatorname{cl}}

\newcommand{\A}{\mathcal{A}}

\newcommand{\K}{\mathcal{K}}
\newcommand{\G}{\mathcal{G}}

\newcommand{\C}{\mathcal{C}}

\newcommand{\R}{\mathbb R}
\newcommand{\Z}{\mathbb Z}
\newcommand{\ep}{\varepsilon}

\theoremstyle{plain} 
\begingroup 

\newtheorem{thm}{Theorem}[section] 
\newtheorem{cor}[thm]{Corollary}
\newtheorem{lem}[thm]{Lemma} 
\newtheorem{prop}[thm]{Proposition}
\newtheorem{defn}[thm]{Definition} 
\newtheorem{rem}[thm]{Remark} 

\endgroup

\theoremstyle{definition}
\theoremstyle{remark}
\numberwithin{equation}{section}

\def\R{{\mathbb R}}
\def\Z{{\mathbb Z}}
\def\C{{\mathbb C}}
\def\D{{\mathbb D}}

\def\Q{{\mathbb Q}}
\def\T{{\mathbb T}}
\def\A{{\mathbb A}}

\begin{document}

\title{Regions of instability for non-twist maps}
\author{John Franks and  Patrice Le Calvez }
\date{}

\begin{abstract}
In this paper we consider an analog of the regions of instability 
for twist maps in the context of area preserving diffeomorphisms 
which are not twist maps.  Several properties analogous to those
of classical regions of instability are proved.
\end{abstract}

\maketitle

\bigskip
\section{Introduction and Notation}

The objective towards which this article is directed is a good
topological picture of the dynamics of $C^r$ generic area preserving
diffeomorphisms of surfaces.  We address the case where the surface is
either the annulus or $S^2$ and $r \ge 3.$ Even in this case our
objective is an ambitious one. Despite some substantial progress, many
fundamental questions remain unanswered.  Perhaps the most important
is the question of whether there is a $C^r$ generic set of area
preserving diffeomorphisms for which the periodic points are dense.
It is not difficult to see that for generic diffeomorphisms this
question is equivalent to the question of whether the union of
all homoclinic points for all hyperbolic periodic orbits forms a 
dense set.

One of our main results, Theorem \ref{dense_stable}, asserts that for
any $r \ge 3$ there is a $C^r$ generic set of diffeomorphisms of $S^2$ with
the property that the union of stable manifolds of all hyperbolic
periodic points is a dense set. Of course, one can simultaneously have the
union of all unstable manifolds of hyperbolic periodic points dense, but
this does not imply the desired result about homoclinic points.

A guiding principle, in pursuing a topological understanding of area
preserving diffeomorphisms has been to establish
analogs of the classical topological results about twist maps (primarily due to 
G.~D.~Birkhoff) in a context where the twist hypothesis is removed,
but generic properties are assumed.  

Our most successful result in this direction is the formulation of a
generalization of Birkhoff's notion of ``region of instability'' and
the proof that many of their key properties remain true for generic
area preserving diffeomorphisms of the annulus or $S^2.$ Our precise
definition of region of instability for an annulus diffeomorphism is
given in Definition \ref{region_d}, but in particular it is an
essential open invariant annulus which contains no essential compact
connected invariant sets which are periodic point free.  In this sense
it is very like the Birkhoff definition for a region of instability
for a twist map which would require that it contain no essential
invariant curves.  For a generic twist map the two definitions
co-incide.

We show that inside any region of instability there is a hyperperbolic
periodic point $z$ with the property that $K(z)$, the closure of
$W^s(z),$ is essential (in the annulus) and contains both components
of the frontier of the region of instability.  There are many such
``essential'' periodic points and we show in Theorem \ref{8.6} and
Corollary \ref{8.7} that $K(z)$ is independent of the essential
periodic point $z$ and that we could have equally well defined it to
be the closure of the unstable manifold $W^u(z).$ Further study of the
dynamics on $K(z)$ is crucial for a good topological understanding of
area preserving diffeomorphisms.  In particular if, in the generic
setting, homoclinic points of sucha point $z$ were always dense in
$K(z)$ this would answer affirmatively the question of whether
periodic points are generically dense for area preserving
diffeomorphisms of the annulus.

One of the most striking results we establish for regions of
instability is a generalization of known results (see Mather \cite{M2})
for connecting orbits of twist maps. In Theorem
\ref{connecting_t} we show that in a region of instability there are
always connecting orbits between the two components of its frontier.
More precisely if the two components of the frontier are $A$ and $B$,
then there is an orbit whose $\alpha$-limit set is in $A$ and whose
$\omega$-limit set is in $B$, as well as one whose $\alpha$-limit set
is in $B$ and whose $\omega$-limit set is in $A$.  While Mather's
proof of the analogous result for twist maps uses variational 
techniques, he has also made seminal contributions to the topological
methods on which our work is based (see \cite{M1}, \cite{M3}).

In the subsequent text, we write $\T^1$ for the one dimensional torus
and $\A=\T^1\times\R$ for the infinite annulus. We note $N$
(resp. $S$) the upper (resp. lower) end of $\A$ and $\widehat\A$ the
end compactification of $\A$ which is a topological sphere, i.e., a
topological set homeomorphic to the two dimensional sphere
\[
S^2=\{(x_1,x_2,x_3)\in\R^3\enskip\vert\enskip x_1^2+x_2^2+x_3^2=1\}.
\]

We consider the covering projection:
\begin{align*}
\pi: \R^2&\to\A\cr
(x,y)&\mapsto(x+\Z,y)
\end{align*}
of the universal covering space $\R^2$ and the Deck
transformation
\begin{align*}
T: \R^2&\to\R^2\cr
(x,y)&\mapsto(x+1,y).
\end{align*}

We give the same name to the first projection
\[
p_1:(x,y)\mapsto x
\]
defined on the plane or the annulus.

A set $K$ of $\A$ is called an {\it essential continuum} if it is
compact, connected, and if its complement has two unbounded
components. We write $U_K^+$ for the upper one and $U_K^-$ for the
lower one. The sets 
\[
\widehat U_K^+=U_K^+\cup \{N\}\enskip\enskip{\rm and}\enskip\enskip 
\widehat U_K^-=U_K^-\cup \{S\}
\]
are simply connected domains of
$\widehat \A$. If $\A \setminus K$ has no other component we will say that $K$ is a
{\it filled essential continuum}.

We say that a homeomorphism of $\A$ preserves area if its preserves
the measure given by the volume form $dx\wedge dy$.

An open set $U$ of $\A$ is called an {\it essential annulus} if it is
homeomorphic to $\A$
and if the inclusion $i: U\hookrightarrow \A$ induces a one to one map
\[
i_*:H_1(U,\R)\to H_1(\A,\R)
\]
on the first homology group. If $U$ is bounded, its frontier has
two connected components, we denote $\partial ^+U$ (resp.
$\partial^-U$) the upper (resp. the lower one).

We will consider 
\[
S^1=\{z\in\C\enskip\vert\enskip \vert z\vert =1\}
\]
and
\[
\D=\{z\in\C\enskip\vert\enskip \vert z\vert <1\}.
\]

\section{Prime end theory}

Let $U$ be a non trivial simply connected domain of a topological
sphere $\Sigma$ (that means that the complement has at least two
points). We can define the prime end compactification of $U$,
introduced by Caratheodory \cite{C}, by adding the circle $S^1$.
A very good exposition of the the theory of prime ends in modern
terminology and in the context we need can be found 
in the paper \cite{M3} of John Mather.  This paper also extends the
theory to surfaces other than the plane or $S^2$ though we shall
not make use of this.

The prime end compactification can be defined purely topologically but
has another significance if we put a complex structure on $\Sigma$. We
can find a conformal map $h$ between $U$ and the open disk $\D$ and we
put on $U\sqcup S^1$ the topology (which, up to homeomorphism of the
resulting space, is independent of $h$) induced from the natural
topology of $\overline
\D$ by the bijection 
\[
\overline h:U\sqcup S^1\to \overline \D
\]
equal to $h$ on $U$ and
to the identity on $S^1$.

If $F$ is a homeomorphism of $\Sigma$ which leaves $U$ invariant, then
the map $F_{\vert U}$ extends to a homeomorphism of the prime end
compactification $U\sqcup S^1$ of $U$. If $F_{\vert U}$ is orientation
preserving then its restriction to $S^1$ preserves orientation and we
can define its rotation number which is an element of $\T^1$. We
denote it by $\rho_U$ and call it the {\it rotation number} of $U$.

Every homeomorphism $F$ of $\A$ homotopic to the identity extends
naturally to an orientation preserving homeomorphism $\widehat F$ of
$\widehat \A$ which fixes the ends $N$ and $S$. If $K$ is an invariant
essential continuum we write $-\rho_K^+$ (resp.  $\rho_K^-$) for the
rotation number of $\widehat U_K^+$ (resp. $\widehat U_K^-$) defined
by $\widehat F$.

Every lift $f$ of $F$ to the plane induces naturally a lift of
$F_{\vert U_K^+}$ to the universal covering space. This lift can be
extended by a lift $\varphi$ of $\Phi$ to the universal covering space
of $S^1$ which admits a real rotation number $\widetilde \rho_K^+$
such that $\widetilde \rho_K^++\Z=\rho_K^+$. In the same way we can
define $\widetilde\rho_K^-$.

By prime ends theory we can compactify any invariant and bounded
essential open annulus, by adding two copies $C^+$ and $C^-$ of $S^1$
corresponding respectively to the upper and to the lower end of
$U$. We will get naturally two rotation numbers $\rho_U^+$ and
$\rho_U^-$ and for a lift $f$ two real rotation numbers
$\widetilde\rho_U^+$ and $\widetilde\rho_U^-$

We will often make use of the following classical result about prime ends:

\begin{prop} \label{fixed_end_prop} Let $F$ be an
orientation preserving homeomorphism of a topological sphere $\Sigma$
and $U$ a non trivial simply connected invariant domain. We suppose
that:

\begin{enumerate}
\item[i)] the rotation number $\rho_U=p/q+\Z$ is rational, and 

\medskip
\item[ii)] there is a neighborhood of $\partial U$ in
$\overline U$
which doesn't contain the positive or the negative orbit of any wandering
open set.
\end{enumerate}
\noindent Then $\partial U$ contains a fixed point of $F^q$.  There is
an analogous result for the case that $U$ is an open annulus.
\end{prop}

\begin{proof} 
Taking $F^q$ instead of $F$, we can
suppose that the rotation set $\rho_U$ is equal to zero.  It follows
that there is a fixed point $z\in S^1$ of the extended homeomorphism
in the prime end compactification of $U$. 

{From} prime end theory we can find a sequence of open arcs
$(\gamma_n)_{n=1}^\infty$ in $U$ of diameter $<1/n$ each of which
extends to a closed arc in the prime end-compactification, with end
points which belong to $S^1$ and which are close to $z$, one on each
side. Moreover, these arcs have the property that the closures
$\overline \gamma_n$ converge to the point $z$ in the prime end
compactification as $n$ goes to infinity.  

Each of these arcs divides $U$ into two simply connected domains.
We denote by $U_n$ the small one (i.e. $z$ belongs to the closure of
$U_n$ in the prime end compactification).  We wish to show that
$F(\gamma_n) \cap \gamma_n \ne \emptyset.$
If the arc $\gamma_n$ doesn't meet its image under $F$, the fact that
$z$ is fixed implies that one of the inclusions $F(U_n)\subset U_n$ or
$F^{-1}(U_n)\subset U_n$ is true. In the first case we can find an
open set $V$ contained in $U_n \setminus F(U_n)$ which is wandering and
has a positive orbit contained in $U_n$.  In the second case we can
construct a wandering domain whose negative orbit is contained in $U_n$.
Since these conditions contradict our hypothesis we conclude that
$F(\gamma_n) \cap \gamma_n \ne \emptyset.$

This fact together with the fact that the diameters of the arcs 
$\gamma_n$ tend to zero implies that any limit point in $\partial U$
is a fixed point of $F$.
\end{proof}

\bigskip
\begin{rem} \label{fixed_end_rem}  We note that in the proof above
we showed that for any fixed prime end $z$ and any sequence of arcs
$\gamma_n$ defining it, we have $F(\gamma_n) \cap \gamma_n \ne \emptyset.$
We shall make use of this fact subsequently.
\end{rem}

We deduce immediately the following corollaries of Proposition
\ref{fixed_end_prop}:

\begin{cor}Let $F$ be an orientation preserving
homeomorphism of a topological sphere
$\Sigma$ and $U$  a non trivial simply connected invariant
domain. We suppose that:

\begin{enumerate}
\item[i)] the rotation number $\rho_U=p/q+\Z$ is rational, and

\item[ii)] there is  no wandering point in $U$.
\end{enumerate}

Then $\partial U$ contains a fixed point of $F^q$.
\end{cor}

\begin{cor}Let $F$ be an area preserving homeomorphism of $\A$ homotopic
to the identity and $K$ an invariant continuum. We suppose that
$\rho_K^+=p/q+\Z$ (resp. $\rho_K^-=p/q+\Z$) is rational. Then $K$
contains a fixed point of $F^q$, more precisely the fixed point is in
$\partial U_K^+$ (resp.  $\partial U_K^-$).
\end{cor}

\bigskip
\section{Rotation number of invariant measures}

Let $F$ be a homeomorphism of $\A$ homotopic to the identity and $f$ a
lift of $F$ to $\R^2$. If $\mu$ is an invariant Borel probability
measure we can define its rotation number
\[
\rho(\mu)=\int_{\A} \phi\,d\mu,
\]
where $\phi:\A\to\R$ is the function lifted by
the function $p_1\circ f-p_1$.

If $O$ is a periodic orbit of $F$ with $q$ elements, there exist an integer
$p\in\Z$ such
that $f^q(z)=T^p(z)$ for every $z\in\pi^{-1}(O)$. We will say that $O$ is a
periodic
orbit of {\it type} $(p,q)$. The number
$p/q$ is called the rotation number of $O$, it is equal to the rotation
number of the
invariant measure
\[
\mu=\frac{1}{q}\sum_{z\in O} \delta_z,
\]
where $\delta_z$ is the
Dirac measure
concentrated on $z$.

\medskip
If $K$ is a non empty compact invariant set, the set of Borel probability
measures whose support is in $K$, is a compact set for the weak
topology and the subset $M(K)$ of the invariant measures is a non empty
convex and closed
set. We can define the {\it rotation set} of $K$
\[
\rho(K)=\left\{\rho(\mu)\enskip\vert\enskip \mu\in M_K\right\},
\]
which is a
non-empty segment $[\alpha,\beta]$ of
$\R$. We have the following classical result:

\begin{prop} For every $\ep > 0$ there exists $N\geq 0$ such that for every
$z\in\pi^{-1}(K)$ and for
every $n\geq N$ we have:
\[
\alpha-\varepsilon<\frac{p_1\circ f^n(z)-p_1(z)}{n}<
\beta+\varepsilon.
\]
\end{prop}

\begin{proof}
Let's suppose for example that we can find a real number
$\varepsilon>0$, a sequence of integers $(n_k)_{k\geq 0}$ converging
to $+\infty$ and a sequence $(z_k)_{k\geq 0}$ in $\pi^{-1}(K)$ such
that 
\[
\frac{p_1\circ f^{n_k}(z_k)-p_1(z_k)}{n_k}\geq \beta+\varepsilon.
\]
We have the inequality
\[
\int_{\A}\phi\,d\mu_k\geq \beta+\varepsilon
\]
where
\[
\mu_k=\frac{1}{n_k}\sum_{i=0}^{n_{k}-1}\delta_{F^{i}(\pi(z_{n_k}))}.
\]

Any closure value $\mu$ of the sequence $(\mu_k)_{k\geq 0}$ belongs to
$M(K)$ and satisfies the inequality

\[
\int_{\A}\phi\,d\mu \geq \beta +\varepsilon.
\]
 We have got a
contradiction.
\end{proof}

We obtain the following corollary which applies for example when $K$ is an
invariant
essential simple closed curve.

\begin{cor}
If $\rho(K)=\{\alpha\}$ the sequence 
\[
\frac{p_1\circ f^n-p_1}{n}
\]

converges uniformly on $\pi^{-1}(K)$ to the constant function
$\alpha$.
\end{cor}

In the same way we can define the rotation number of an invariant
probability measure when $F$ is a homeomorphism of the closed annulus
$\T^1\times [0,1]$ homotopic to the identity, for a lift $f$ to
$\R\times [0,1]$. The rotation set of $f$ is the set of the rotation
numbers of all the invariant measures. It is a segment which contains
the rotation number of the two boundary curves of the annulus. We will
use the following result proved in
\cite{F1}:

\begin{thm} \label{birkhoff_thm} 
Let $F$ be an area preserving homeomorphism of the closed annulus
$\T^1\times [0,1]$ homotopic to the identity and $f$ a lift of $F$ to
$\R\times [0,1]$.  For every rational number $\rho $ contained in the
rotation set of $f$, there exists a periodic orbit of type $(p,q)$,
where $\rho=p/q$ and $p$ and $q$ are relatively prime.
\end{thm}

\bigskip
\section{Indices}

Let $F$ be an orientation preserving homeomorphism of a topological
sphere $\Sigma$ and $D$ a closed disk of $\Sigma$ (i.e. a set
homeomorphic to $\overline \D$).  If there is no fixed point of $F$ in
$\partial D$, we can define the index $i(F,D)$ of $D$ (see \cite{D}
for example). This index is equal to the sum of the Lefschetz indices
$i(F,z)$ of the fixed points $z$ contained in $D$ if there are a
finite number of fixed points.

If we can find a finite family of disjoint closed disks $(D_k)_{1\leq
k\leq n}$ such that every fixed point of $F$ is contained in the
interior of one these disks, we have the Lefschetz formula:

\[
\sum_{k=1}^n i(F,D_k)=2
\]

If $U$ is a non trivial simply connected domain of $\Sigma$ with a non
zero rotation number, we can find a closed disk $D\subset U$ whose
interior contains all the fixed points in $U$ and such that
$i(F,D)=1$. We deduce that $U$ contains at least one fixed point.

\bigskip
\section{Invariant essential continua and periodic orbits.}

We begin with a useful extension lemma

\begin{lem}\label{extension_lem}
Let $F$ be a homeomorphism of $\A$ homotopic to the identity and $f$ a
lift of $F$ to $\R^2$. We suppose that $K$ is an invariant essential
continuum and that $\widetilde\rho_K^+>0$.  Then there exists a
homeomorphism $F'$ of $\A$ homotopic to the identity and an essential
simple closed curve $C^+\subset U_K^+$ such that:

\begin{enumerate}
\item[i)]
the map $F'$ coincides with $F$ on the complement of a closed set contained
in $U_K^+$;
\item[ii)]
the curve $C^+$ is invariant by $F'$ and its rotation number is $>0$
(for the lift $f'$ which coincides with $f$ on $\pi^{-1}(K)$);
\item[iii)]
the invariant set $U_K^+\cap \overline {U_{C^+}^-}$ of
$F'$ contains no fixed points of rotation number $0$;
\item[iv)]
the map $F'$ is area-preserving if $F$ is.
\end{enumerate}
\end{lem}

\begin{proof} We can find a complex structure on
$\widehat \A$ which
defines the same differential structure on $\A$. Let's consider a conformal map
$h$ between the simply connected domain
$\widehat U_K^+$ and $\D$ and the maps
\begin{align*}
r:U_K^+&\to ]0,1[\cr
 z&\mapsto \vert h(z)\vert\cr
\end{align*}

and
\begin{align*}
\theta: U_K^+&\to S^1\cr
z&\mapsto \frac{h(z)}{\vert h(z)\vert}\cr
\end{align*}

We consider also a lift

\[
\widetilde \theta:\pi^{-1}(U_K^+)\to \R
\]
of $\theta$ to the universal covering space and the map
\[
\widetilde r: \pi^{-1}(U_K^+) \to ]0,1[
\]
given by $\widetilde r(z) = r(\pi(z)).$

For every $t\in]0,1[$ we consider the essential simple closed curve
$C_t=r^{-1}(\{t\})$ of
$U_K^+$. Let's fix $t_0<1/2$ such that the curve $C^+=C_{t_0}$ is above the
image of
$C_{1/2}$. We can find a homeomorphism between the closed annulus $A_0$
delimited by
$C_{1/2}$ and
$C_{t_0}$ and the closed annulus $A_1$ delimited by $F(C_{1/2})$ and
$C_{t_0}$, which is equal to $F$ on $C_{1/2}$ and to the identity on
$C_{t_0}$. We denote
by $F_1$ its extension to the infinite annulus which is equal to $F$ below
$C_{1/2}$ and
to the identity above $C_{t_0}$. If $F$ is area preserving the fact that
$K$ is invariant
tells us that the two annuli $A_0$ and $A_1$ have the same area and we can
choose the
homeomorphism $F_1$ to be area-preserving (so is the map $F_1$).

We write $f_1$ for the lift of $F_1$ which coincides with $f$ on
$\pi^{-1}(K)$. Using the condition
$\widetilde\rho_K^+>0$, we can find $\varepsilon\in]0,1/4[ $ such that
\[
\widetilde r(z)>1-2\varepsilon\enskip\Rightarrow\enskip\widetilde\theta\circ
f_1(z)>\widetilde \theta (z).
\]

Using an argument of compactness, we can find  $M>0$
such that
\[
\widetilde r(z)\geq t_0\enskip\Rightarrow\enskip\widetilde\theta\circ
f_1(z)\leq
\widetilde\theta (z)+M.
\]

Consider a $C^{\infty}$ decreasing map 
\[
\lambda: ]0,1[ \to [0,+\infty[
\]
 
equal to $0$ on
$[1-\varepsilon,1[$ and let's define the map $F_2$ equal outside of
$U_K^+$ to the identity and on $U_K^+$ to the time one map
of the Hamiltonian flow of
$\lambda\circ r$ for the usual symplectic structure on
$\A$. That means the flow defined by the equations:

\begin{align*}
x'&=\lambda'(r(x,y))\frac{\partial r}{\partial y}(x,y)\\
y'&=-\lambda'(r(x,y))\frac{\partial r}{\partial y}(x,y).
\end{align*}

The map $F'=F_1\circ F_2$
satisfies the conditions {\it i)} and  {\it iv)} of the lemma. It also
fixes the curve
$C^+$. We can choose the function
$\lambda$ such that its satisfies the condition {\it iii)} and the second
part of the
condition {\it ii)}.

\medskip

Indeed, if $f_2$ is the lift of $F_2$ equal to identity on $\pi^{-1}(K)$,
we have
\[
\widetilde\theta\circ f_2(z)\geq \widetilde\theta(z)\enskip\enskip{\rm
and}\enskip\enskip \widetilde r\circ f_2(z)=\widetilde r(z)
\]
 for every
$z\in\pi^{-1}(U_K^+)$, so we have the relation
\[
\widetilde r(z)>1-2\varepsilon\enskip\Rightarrow\enskip\widetilde\theta\circ
f'(z)>\widetilde \theta(z).
\]

If we choose the function $\vert\lambda'\vert$ sufficiently large on $[r_0,
1-2\varepsilon]$ we will get
\[
\widetilde r(z)\in[r_0,1-2\varepsilon]\enskip\Rightarrow\enskip\widetilde\theta\circ
f'(z)>\widetilde \theta (z).
\]

\end{proof}

\bigskip
We have a similar result in the cases where $\widetilde\rho_K^+<0$,
$\widetilde\rho_K^->0$ and
$\widetilde\rho_K^-<0$.

\medskip
 The main result of this paragraph is the following.

\begin{prop} \label{filled_prop}
Let $F$ be an area preserving homeomorphism of $\A$ homotopic to the
identity and $f$ a lift of $F$ to $\R^2$. If $K$ is an invariant
essential continuum which doesn't contain any periodic orbit of $F$,
then:

\begin{enumerate}
\item[i)] the set $K$ is a filled essential continuum;

\item[ii)] there exists $\alpha\in\R \setminus \Q$ such that
$\rho(K)=\{\alpha\}$;

\item[iii)] $\widetilde\rho_K^+=\widetilde\rho_K^-=\alpha$.
\end{enumerate}
\end{prop}

\begin{proof}
We first prove assertion {\it i)} using the results of the paragraphs
2 and 4.  Suppose that $\A \setminus K$ has a bounded component
$U_0$. It is periodic and we can suppose it fixed by taking an iterate of
$F$ instead of $F$. Consider the extension $\widehat F$ of $F$ to the
end compactification $\widehat \A$. 

If $U$ is any invariant bounded component of $\A \setminus K$ then it
is topologically an open disk and must contain a fixed point since $F$
is area preserving.  The fact that $K$ doesn't contain any fixed
points implies that the set of invariant bounded components of $\A
\setminus K$ is finite, since otherwise there would be a sequence of
fixed points with limit in $K$. Moreover we can find a finite family
$(D_k)_{1\leq k\leq n}$ of closed disks, one in each fixed component
of $\widehat \A \setminus K$, such that every fixed point of $\widehat
F$ is contained in the interior of one these disks and such that
$i(\widehat F,D_k)=1$ for every $k$. There are at least three
components: $\widehat U_K^+$, $\widehat U_K^-$ and $U_0.$ This
contradicts the Lefschetz formula.

\bigskip
We next prove assertions {\it ii)} and {\it iii)}. Using Proposition
\ref{fixed_end_prop}, we know that the rotation numbers $\widetilde\rho_K^+$ and
$\widetilde\rho_K^-$ are irrational. Let's consider a rational number
$p/q$. Applying the extension lemma to $F^q$ and to the lift $f^q\circ
T^{-p}$, one can construct an area preserving homeomorphism $F'$ of
$\A$ which coincides with $F^q$ on a neighborhood of $K$ and which
fixes two essential simple curves $C^+\subset U_K^+$ and $C^-\subset
U_K^-$ such that:

\medskip
\noindent
$\bullet$\enskip\enskip there is no fixed point of rotation number $0$
\ in $U_K^+\cap\overline {U_{C^+}^-}$;

\medskip
\noindent
$\bullet$\enskip\enskip there is no fixed point of rotation number $0$
\ in $U_K^-\cap\overline {U_{C^-}^+}$;

\medskip
\noindent$\bullet$\enskip\enskip the rotation number of $C^+$ and
$q\widetilde\rho_K^+-p$
have the same sign;

\medskip
\noindent$\bullet$\enskip\enskip the rotation number of $C^-$ and
$q\widetilde\rho_K^--p$
have the same sign.

\bigskip
Using the fact that $K$ is filled and doesn't contain any periodic
orbit we deduce that the map $F'$ has no fixed point of rotation
number $0$ in the closed annulus defined by the curves $C^+$ and
$C^-$. {From} Theorem \ref{birkhoff_thm} we deduce that:

\medskip
\noindent$\bullet$ \enskip\enskip the rotation set $\rho(K)$ doesn't
contain $p/q$;

\medskip
\noindent$\bullet$ \enskip\enskip the numbers $\widetilde\rho_K^-$ and
$\widetilde\rho_K^+$ are not separated by
$p/q$.

The conclusion follows immediately.
\end{proof}

\begin{defn} \label{irr_cont_defn} Such a set will be called an 
{\it irrational invariant essential continuum} and $\alpha=\rho(K)$
its rotation number.
\end{defn}

\medskip
In the same way, we can prove the following result.

\begin{prop} \label{per_pt_prop}Let $F$ be an area
preserving homeomorphism of $\A$ homotopic to the identity, $f$ a lift
of $F$ to $\R^2$ and $K$ an invariant essential continuum. For every
rational number $\rho $ contained in $\rho(K)$, there exists a
periodic orbit of type $(p,q)$ contained in $K$, where $\rho=p/q$ and
$p$ and $q$ are relatively prime.
\end{prop}

\begin{proof} We give a proof by contradicting
the assumption that $F$ has no periodic points of type $(p,q)$ in $K$.
If this is assumed then $F$ has no periodic points of type $(p,q)$ in 
some neighborhood $V_0$ of $K$ since the periodic points of type $(p,q)$ 
form a closed set.  It then follows from Proposition \ref{fixed_end_prop} that
$\widetilde\rho_K^+ \ne p/q$ and $\widetilde\rho_K^- \ne p/q$.

Using Lemma \ref{extension_lem} we construct an area preserving map
$F'$, which coincides with $F^q$ on some neighborhood $V \subset V_0$
of $K$, which leaves invariant two essential simple closed curves
$C^+\subset U_K^+$ and $C^-\subset U_K^-$, and which has no fixed
points of rotation $0$ in $\pi^{-1}(V)$ for the lift equal to
$f^q\circ T^{-p}$ on $\pi^{-1}(K)$. By Theorem \ref{birkhoff_thm},
there exists a fixed point of $F'$ of rotation number $0$ in the
annulus $B$ delimited by $C^+$ and $C^-$. Since it is not in $K$ and
not in $V$ this fixed point must lie in a bounded invariant component
of $B \setminus K$ and hence a bounded invariant component of $A
\setminus K.$

Let $p: B_0 \to B$ be a finite covering map of degree greater than the
rotation number of any periodic point of $F'$.  Then there is a lift
$G : B_0 \to B_0$ of $F'$ such that $z \in B_0$ is a fixed point of
$G$ if and only if $p(z)$ is a fixed point of $F'$ of rotation number
zero.  Let $K_0 = p^{-1}(K)$ so that $G$ has no fixed points in $K_0.$
Since any component of $B \setminus K$ containing a fixed point of $F'$ of
rotation number zero is an invariant open topological disk it follows that any
component of $B_0  \setminus  K_0$ containing a fixed point of $G$ is also an
invariant open topological disk.

As in the proof of Proposition \ref{filled_prop}, the fact that $K_0$ and the
boundary of $B_0$ don't contain any fixed points of $G$ implies that the set
of bounded components of $B_0 \setminus K_0$ containing fixed points is
finite.  Moreover we can find a finite family $(D_k)_{1\leq k\leq n}$
of closed disks, one in each component of $B_0 \setminus K_0$ containing such a
fixed point, such that every fixed point of $G$
is contained in the interior of one these disks and such that
$i(G, D_k)=1$ for every $k$. There is at least one such disk since 
there was at least one fixed point of $F'$ of rotation number zero.
This again contradicts the Lefschetz formula for $G: B_0 \to B_0.$
\end{proof}

\bigskip
\section{Generic properties of diffeomorphisms of the sphere.}

We consider a smooth surface $\Sigma$ diffeomorphic to $S^2$ and a
volume form $\omega$ on it. For every integer $k\geq 1$ we write ${\rm
Diff}^r_{\omega}(\Sigma)$ for the set of $C^r$ diffeomorphisms which
preserves $\omega$. We put on that set the natural $C^r$-topology. If
$z$ is a periodic point of period $q$ of $F\in{\rm
Diff}_{\omega}^r(\Sigma)$, we have three possibilities:

\begin{enumerate}
\item[$\bullet$] the eigenvalues of $DF^q(z)$ are both equal to $1$ or
$-1$: the point
$z$ is {\it degenerate};

\item[$\bullet$] the eigenvalues of $DF(z)$ are two different
conjugate complex numbers of modulus $1$: the point $z$ is {\it elliptic};

\item[$\bullet$] the eigenvalues $\lambda$ and $\mu$ of
$DF(z)$ are real and satisfy $\vert\lambda\vert<1<\vert \mu\vert$: the point  $z$ is {\it hyperbolic} or a {\it saddle point}.
\end{enumerate}

The {\it stable (resp. unstable) branches} of a hyperbolic periodic
point $z$ are the two components of $W^s(z) \setminus \{z\}$
(resp. $W^u(z) \setminus \{z\}$) for the topology induced by a parameterization
of $W^s(z)$ (resp. $W^u(z)$). Let's recall that a {\it saddle
connection} is a branch $\Gamma$ which is a stable branch of a
hyperbolic periodic point $z$ and an unstable branch of a hyperbolic
periodic point $z'$ (perhaps equal to $z$).

\begin{lem}\label{branch_lem}
Suppose $F$ is an area preserving diffeomorphism of a connected open
subset $U$ of $S^2$ and $K \subset U$ is a compact connected
$F$-invariant set.  If $\Gamma$ is a branch of a hyperbolic periodic
point of $F$ then either $\Gamma \subset K$ or $\Gamma$ is entirely
contained in a single component of $U \setminus K.$
\end{lem}

\begin{proof}
Suppose $\Gamma$ has non-empty intersection with both $K$ and $U
\setminus K.$ Let $U_0$ be a component of $U \setminus K$ which
contains a point of $\Gamma$.  Then $U_0$ is a periodic domain since
every component of the complement of $U \setminus K$ is.  There
is an an open arc $\alpha \subset \Gamma \cap U_0$ whose endpoints are
in $K$.

We note that if $i>0$, then $F^i(\alpha)$ and $\alpha$ are disjoint
since $F^i(\alpha)$ cannot be a subset of $\alpha$ (as $\alpha$ is in
a branch) and $F^i(\alpha)$ cannot contain an endpoint of $\alpha$.

Let $V$ be the component of $S^2 \setminus K$ which contains $U_0.$
Since $V$ is simply connected $\alpha$ separates it into two non-empty
simply connected domains $V_1$ and $V_2$ each of which has non-empty
intersection with $U_0$.  It is not difficult to show that this
implies $\alpha$ separates $U_0$ into two connected open sets $U_1 =
U_0 \cap V_1$ and $U_2 = U_0 \cap V_2$.  

Choose $k>0$ so that $F^k(U_0)=U_0$. The map $F$ is area preserving, so
$U_1\times U_2$ is a non wandering open set for $F^k\times F^k$. We deduce
that there exists an integer $n\geq 1$ such that $F^{nk}(U_1)\cap
U_1 \ne \emptyset$ and $F^{nk}(U_2)\cap U_2 \ne \emptyset$. Moreover the
closure in $U_0$ of  $F^{nk}(U_1)$ is not contained in $U_1$, it meets
the boundary $\alpha$. Similarly $\alpha$ meets the closure of
$F^{nk}(U_2)$ in $U_0$. We deduce that $\alpha $ meets its image by
$F^{nk}$, which is impossible.
We conclude it is not possible for
$\Gamma$ to have non-empty intersection with both $K$ and its
complement.  If $\Gamma$ lies in the complement of $K$ it must be in a
single component because it is connected.
\end{proof}

We will say that a periodic point $z$ of period $q$ is {\it Moser
stable } if it admits a fundamental system of neighborhoods which are
closed disks $D$ such that $F_{\partial D}$ is minimal.

We next present the following result of Mather \cite{M1}.

\begin{thm} \label{mather_thm}(Mather)
Let $F\in{\rm Diff}_{\omega}^r(\Sigma)$ such that:
\begin{enumerate}
\item[$\bullet$] there is no degenerate periodic point,

\item[$\bullet$] every elliptic periodic point is Moser
stable, and

\item[$\bullet$] there are no saddle connections.
\end{enumerate}

\noindent Then, we know that:

\begin{enumerate}
\item[i)] the rotation number of any non trivial
simply connected periodic domain is irrational;

\item[ii)] the branches of a hyperbolic periodic
point have the same closure.
\end{enumerate}
\end{thm}

\begin{proof}
We first will prove assertion {\it i)} by contradiction.  If assertion
{\it i)} is false then by replacing $F$ with an iterate we can assume
that there is an invariant simply connected domain $U$ with rotation number zero.
This implies that the extension of $F$ to the prime end compactification of $U$
has a fixed point $z$ on the boundary.  

{From} the theory of prime ends and the proof of Proposition
\ref{fixed_end_prop} we know that there is a sequence of open arcs
$(\gamma_n)_{n=1}^\infty$ in $U$, each of which is the interior of a
closed arc with endpoints in $\partial U$, which define the prime end
$z$.  Moreover, we can assume the sequence of diameters of
$(\gamma_n)$ tends to zero as $n$ tends to infinity, so by choosing a
subsequence we can assume the closures of the arcs converge to a
point $p$ in $\partial U.$  We observed (see Remark
\ref{fixed_end_rem}) that $F(\gamma_n) \cap \gamma_n \ne \emptyset$
for each $n$ and hence $p$ is a fixed point of $F$.

We note first that this fixed point $p$ cannot be a Moser stable
point.  This is because a Moser stable point is in the interior of an
arbitrarily small $F$ invariant disk on the boundary of which $F$ is
minimal.  Thus the boundary of such a disk must intersect the
complement of $U$ and hence, by minimality, the boundary is disjoint from $U$.
This would imply that the entire disk is disjoint from $U$ which is
impossible.

It follows that the point $p$ must be a hyperbolic fixed point.  By
the Hartman-Grobman theorem there is a neighborhood $V$ of $p$ on which
there are topological co-ordinates $(x,y)$ with $p$ at the origin and
in which $F(x,y) = (2x, y/2)$.  We will use these local coordinates on
$V$ and refer to the four {\it local branches} at $p$ by which we mean
the intersection of $V$ with the positive and negative $x$ and $y$
axes.

The contradiction we wish to derive will follow from the following lemma

\begin{lem}
Replacing the sequence $\{\gamma_n\}$ with a subsequence, we may conclude
that there is a neighborhood $W \subset V$ of $p$ such that 
one component of $U \cap W$ contains the arcs ${\gamma_n}$ for 
$n$ sufficiently large, and this component consists of one of 
the following:

\begin{enumerate}
\item[a)] the intersection of $W$ with a single open quadrant in
the $(x,y)$ co-ordinates, or

\item[b)] the intersection of $W$ with two adjacent open quadrants 
and the local branch between them, or

\item[c)] the intersection of $W$ with three open quadrants 
and the two local branches between them, or

\item[d)] the complement in $W$ of $p$ and a single local branch.
\end{enumerate}
\end{lem}

\begin{proof}
We first observe that if $\gamma$ is a local branch of $p$ in $V$ then
by Lemma \ref{branch_lem} $\gamma$ must be entirely contained in $U$ or
disjoint from $U$.  If it is contained in $U$ there is an $x \in \gamma$ 
such that a neighborhood of the interval $[x, F(x)]$ is in $U$.
The collection of all iterates of this neighborhood contain the
intersection of a smaller neighborhood $W_0$ with the two open quadrants in 
$(x,y)$ co-ordinates of which $\gamma$ is part of the boundary.
Thus if one point of a local branch is in $U$ then, perhaps in a smaller
neighborhood, the entire branch and the open quadrants on either side of
it are in $U$.  We shall make repeated use of this property.

The first use of this fact is to observe that it is not possible for
all four local branches at $p$ to be in $U$ as this would imply that
there is a neighborhood $W$ of $p$ such that $W \setminus \{p\}
\subset U.$ This would mean that $p \in U$ which is not the case.

We again use the fact that if a local branch intersects $U$ then
(perhaps in a smaller neighborhood) the entire branch and the open
quadrants on either side of it are in $U$.  This implies that either
one of properties {\it b) -- d)} above holds or there is a single
quadrant such that the interiors of infinitely many of the arcs
$\gamma_n$ are in that quadrant, but the two branches which bound this
quadrant are not in $U$.  We complete the proof of the lemma by
showing that this last alternative implies that {\it a)} holds.

We can suppose, without loss of generality, that it is the open first quadrant
\[
Q = \{(x,y) | x>0,\ y>0\}
\]
which contains the interiors of infinitely many of the $\gamma_n$ and
the positive $x$ and $y$ axes are in the complement of $U$.  By
choosing a subsequence we may assume every $\gamma_n$ is in the first
quadrant. We need only show that there is a neighborhood $W \subset V$
of $p$ such that $W \cap Q$ is entirely in $U$.

To do this, for each small $\ep>0$ we consider the arc
$\Gamma_\ep$ made up of the two line segments $\{(x,y) = (t,\ep)
| 0 \le t \le \ep\}$ and $\{(x,y) = (\ep, t) | 0 \le t \le
\ep\}$.  These two segments form two sides of a rectangle $R_\ep$
the other sides of which lie in the $x$ and $y$ axes.  This rectangle
contains infinitely many of the arcs $\gamma_n$ and in particular
points of $U$. Also if $\ep$ is sufficiently small we may assume
this rectangle is disjoint from $\gamma_1$ and in particular cannot
contain all of $U$.  It follows that the arc $\Gamma_\ep$ has
non-empty intersection with $U$.

We want next to show there is an open sub-arc $\beta$ of $\Gamma_\ep$
with the property that it is in $U$ but its endpoints are in
the complement of $U$, and that $F(\beta) \cap \beta \ne \emptyset$.
To see this we note that $\Gamma_\ep$ separates $\gamma_1$ and
$\gamma_N$ for a large $N$.  Hence some component $\beta$ of 
$\Gamma_\ep \cap U$ must separate $\gamma_1$ and
$\gamma_N$ in $U$.  Clearly $\beta$ is an open sub-arc of $\Gamma_\ep$
with endpoints in the complement of $U$.  

Let $V_1$ and $V_2$ be the components of $U \setminus \beta$ which
contain $\gamma_1$ and $\gamma_N$ respectively.  The fact that
$F(\gamma_1) \cap \gamma_1 \ne \emptyset$ and $F(\gamma_N) \cap
\gamma_N \ne \emptyset$ implies $F(V_i) \cap V_i \ne \emptyset$
for $i = 1,2.$ Since $F$ is area preserving $F(V_i)$ cannot be a
proper subset of $V_i$.  We conclude that $F(\beta) \cap \beta \ne
\emptyset$.

Since $\beta \subset \Gamma_\ep$ and $F(\Gamma_\ep) \cap 
\Gamma_\ep$ consists of the single point $(x,y) = (\ep,\ep/2)$
we conclude that both $(\ep,\ep/2)$ and its inverse image 
$(\ep/2,\ep)$ are in $\beta.$  Thus the arc $J_\ep$
made of the two line segments
$\{(x,y) = (t,\ep)
| \ep/2 \le t \le \ep\}$ and $\{(x,y) = (\ep, t) | 
\ep/2 \le t \le \ep\}$ is entirely in $U$.  Similarly
for any $\delta <\ep$ we have $J_\delta \subset U$ and
indeed for any $n \in \Z, \ F^n(J_\delta) \subset U.$ 
This easily implies that $R_\ep \cap Q \subset U$, which means
that property {\it a)} holds.  
\end{proof}

We next observe that this lemma and the assumption that {\it i)} is
false do lead to a contradiction.  This is because from the lemma
we see that there is a neighborhood of the prime end $z$ corresponding
to $p$ in the prime ends of $U$ which consists of accessible prime
ends corresponding to points in either stable or unstable branches of
$z$.  It also follows that there are finitely many fixed prime ends and
they all have this property.  We can conclude that any non-fixed prime
end $e$, which is between the fixed prime ends $z_1$ and $z_2$ 
must correspond to a point on a branch of $p_1$ and a branch of $p_2$
where $p_i$ is the fixed point in $\partial U$ corresponding to 
the prime end $z_i$.  Since this is true for the entire interval
of prime ends between $z_1$ and $z_2$ there is a saddle connection
between $p_1$ and $p_2$ which contradicts our hypothesis.  We conclude
that {\it i)} of the theorem must hold.

Property {\it ii)} of the theorem is a consequence of {\it i)} and
Lemma \ref{branch_lem}.  As before, by taking an iterate we may assume
the point in question, $p$, is fixed.  The closure of one branch of
this point is a compact connected invariant set $K$.  {From} Lemma
\ref{branch_lem} we know that another branch of $p$ is either a subset
of $K$ or disjoint from it.  It cannot be disjoint since if $U$
denotes the component of the complement of $K$ which contains $\gamma$
then $U$ would be an invariant simply connected domain and the fixed
point $p$ would be an accessible point of $\partial U$.  This would mean
that the rotation number of $U$ is rational and hence property {\it i)}
would be contradicted.

We conclude $\gamma \subset K$.  Since we have shown any branch of $p$
is a subset of the closure of any other, property {\it ii)} follows.
\end{proof}

\bigskip
We remark that for $r\geq 3$ the set of diffeomorphisms satisfying the
hypotheses of Theorem \ref{mather_thm} is generic in ${\rm
Diff}^r_{\omega}(\Sigma).$

We now recall a result of Robinson \cite{R2} and Pixton \cite{P}.  It
was first proved by Robinson in the case of fixed points and generalized
by Pixton to periodic points.

\begin{thm} \label{pixton_thm}
For every $r\geq 1$, there exists a generic set $G_r$ in
${\rm Diff}^r_{\omega}(\Sigma)$ such that for every $F$ in $G_r$ the following
properties are satisfied:

\begin{enumerate}
\item[i)] if $z$ is a hyperbolic periodic point, any
stable branch $\Gamma$ of $z$ has non-empty transversal intersection
with any unstable branch $\Gamma'$ of $z$.  Moreover
$\overline\Gamma=\overline {\Gamma'}$.

\item[ii)] if $z$ and $z'$ are two distinct hyperbolic periodic
points such that $z'$ belongs to the closure of a branch of $z$, then
any stable (resp. unstable) branch $\Gamma$ of $z$ meets transversely
any unstable (resp. stable) branch $\Gamma '$ of $z'$.  Moreover
$\overline \Gamma=\overline {\Gamma'}$.
\end{enumerate}
\end{thm}

\bigskip
We make some comments on these theorems.

\bigskip
\noindent $\bullet$
Every periodic connected and compact set which contains a hyperbolic
periodic point, and which is not a single point, contains all the
branches of that point.

\noindent $\bullet$ If we use Mather's theorem with the stronger (and
generic)  hypothesis that any intersections of branches are transverse,
then assertion {\it ii)} of Mather's theorem can be used to derive
assertion {\it i)} of Pixton's theorem (see also Oliveira \cite{O} or Arnaud
\cite{A}). Moreover we have the following easy fact:

\noindent $\bullet$ Any hyperbolic periodic point has a
fundamental system of neighborhoods which are closed disks whose
boundary is contained in $W^s(z)\cup W ^u(z)$.
This property and the classical $\lambda$-lemma (see Theorem 11.1 of
Robinson \cite{R1}) implies the condition {\it ii)} of
Pixton's theorem.

\bigskip
\section{ Essential and inessential hyperbolic periodic points.}

In this section we consider an area preserving $C^1$-diffeomorphism of
$S^2$ such that:

\begin{enumerate}

\item[i)] there are no degenerate periodic points,

\item[ii)] every elliptic periodic point is Moser
stable, and

\item[iii)] the intersection of any branches of hyperbolic periodic points
are transverse and any stable and any unstable branch of the same
periodic point have non-empty intersection.
\end{enumerate}

If $\A_0$ is a connected open subset of $S^2$ whose complement has
exactly two components, neither of which is a single point, we will
call it an {\it annular domain.}  It is homeomorphic to an open
annulus.  If $\A_0$ is an $F$-invariant annular domain then we make
speak of an {\it essential} continuum $K \subset \A_0$ meaning one
which separates the two components of the complement of $\A_0.$

Mather's theorem and Pixton's theorem immediately imply the following corollaries.

\begin{prop}\label{generic_prop}
For a diffeomorphism $F$ of $S^2$ satisfying i) - iii) 
above  we have the following properties:

\begin{enumerate}
\item every periodic simply connected domain has an
irrational prime end rotation number;

\item every periodic annular domain has both prime end
rotation numbers irrational;

\item if $\A_0$ is an invariant annular domain and 
$K \subset \A_0$ is an invariant essential continuum, the rotation 
numbers $\rho_K^+$ and $\rho_K^-$ are irrational;

\item  
every stable branch $\Gamma^s$ and every unstable branch $\Gamma^u$ of a
hyperbolic periodic point have a non-empty transversal intersection
and we have $\overline{\Gamma^s}=\overline{\Gamma^u}$;

\item
if $z$ is a hyperbolic periodic point and $\Gamma^s$ and $\Gamma^u$
are a stable and an unstable branch of $z$ then there is a fundamental
system of neighborhoods of $z$ which are closed disks whose boundary
is contained $\Gamma^s \cup \Gamma^u$;

\item a periodic compact connected set,
which is not a single point, and which contains a hyperbolic periodic point,
contains the branches of that point;

\item if the closure of a  branch of
a hyperbolic periodic point $z$ contains a hyperbolic periodic point
$z'$, then every stable (resp. unstable) branch $\Gamma$ of $z$ meets
transversely any unstable (resp.  stable) branch $\Gamma'$ of $z'$ and
we have $\overline \Gamma=\overline{\Gamma'}$.
\end{enumerate}
\end{prop}
\bigskip

\bigskip
Let $F: \A_0 \to \A_0$ be an area preserving diffeomorphism of an annular
domain which extends continuously to $\cl \A_0 \subset S^2$.  We
consider a hyperbolic periodic point $z$ for $F$.  We will assume the
generic property that each stable branch of $z$ has non-empty
intersection with each unstable branch and all intersections are
transverse.

\begin{defn}
If there is an essential closed curve in $\A_0$ consisting
of a finite number of segments of branches of $z$ we will say that $z$
is {\it essential}, otherwise we say that $z$ is {\it inessential}.
\end{defn}

It is immediate that a hyperbolic periodic orbit $O$ consists
of points which are all essential or all inessential.

\medskip
\begin{prop}\label{inv_domain_prop}
Suppose $z \in \A_0$ is a hyperbolic fixed point of $F$.  Let $U(z)$
be the union of all closed topological disks in $\A_0$ whose
boundaries consist of finitely many segments of the branches of $z$, then

\begin{enumerate}
\item[i)] $U(z)$ is an open, connected, and $F$-invariant set which contains $z$ and
its branches.

\item[ii)]
If $z$ is inessential then $U(z)$ is simply connected.

\item[iii)]
If $z$ is essential then $U(z)$ is annular and essential in $\A_0.$

\item[iv)]
If $V \subset \A_0$ is any other open connected $F$-invariant set which
is either simply connected or annular and essential then
either $U(z) \subset V,\ V \subset U(z),$ or $U(z) \cap V = \emptyset.$
\end{enumerate}
\end{prop}

\begin{proof}
Let $z$ be a hyperbolic fixed point of $F$ and let $U(z)$ be the union
of all closed topological disks in $\A_0$ whose boundaries consist of
finitely many segments of the branches of $z$.  Clearly $F(U(z)) =
U(z).$ Property 3 of Proposition \ref{generic_prop} implies $z$ is in
the interior of $U(z)$.  Since any closed segment on a branch of $z$
can be mapped by a power of $F$ into a small neighborhood of $z$ we
see that these segments are in the interior of one of the closed
disks whose union forms $U(z).$ It follows that any of these closed
disks is contained in the union of the interiors of finitely many
other such disks, so $U(z)$ is also equal to the the union of the
interiors of these disks and is open.  This proves {\it i)}. 

We next show {\it ii)}.  This is equivalent to showing that $z$
inessential implies the complement of $U(z)$ in $S^2$ is connected.
If this is not the case then there is a simple closed curve $\alpha$
in $U(z)$ which separates two components, call them $K_1$ and $K_2$,
of the complement of $U(z).$ Since $U(z)$ is the union of the
interiors of closed disks whose boundaries are branch segments, the
curve $\alpha$ is covered by the interiors of finitely many such
disks.  Let $V$ denote the union of these open disks and let $\beta$
be the frontier of the component of $S^2 \setminus V$ containing
$K_1.$ Then $\beta$ is a simple closed curve which separates $K_1$ and
$K_2$ and consists of finitely many segments of branches of $z$.  If
$\beta$ were essential in $\A_0$ then $z$ would be essential.  Hence
we know $\beta$ bounds a disk in $\A_0$ which must contain one of
$K_1$ or $K_2$.  This would imply that one of $K_1$ or $K_2$ is in
$U(z)$ which is a contradiction.  Hence we have shown that $z$
inessential implies that $U(z)$ is simply connected.

If $z$ is essential then a similar argument shows that the complement
of $U(z)$ has exactly two components so {\it iii)} follows.  More
precisely, the fact that $z$ is essential implies $U(z)$ is not simply
connected so $S^2 \setminus U(z)$ has at least two components and in
particular there are distinct components $K_1$ and $K_2$ each of which
contains one of the components of $S^2 \setminus \A_0$.  If {\it iii)}
does not hold then the complement of $S^2 \setminus U(z)$ has at least
three connected components: $K_1$ and $K_2$ and one more $K_3$ which
is a subset of $\A_0.$ There is a simple closed curve $\alpha$ in
$U(z)$ which separates $K_3$ from $K_1 \cup K_2$ and hence from $S^2
\setminus \A_0$.  As above we can find $\beta_0$ a simple closed curve which
separates $K_3$ and $S^2 \setminus \A_0$ and which consists of
finitely many segments of branches of $z$.  Since $\beta_0$ does not
separate the two components of $S^2 \setminus \A_0$ it bounds a closed
disk in $\A_0$.  The interior of this disk contains $K_3$ but is also
in $U(z)$, a contradiction.  We conclude that $S^2 \setminus U(z)$ has
only two components.

Finally we show {\it iv)}.  If $U(z)$ does not intersect both $V$ and
the complement of $V$ then either $U(z) \subset V$ or $U(z) \cap V =
\emptyset.$ But if $U(z)$ does intersect both $V$ and its complement
then either $V$ is entirely contained in the interior of a disk whose
boundary consists of branch segments (in which case $V \subset U(z)$ )
or $V$ contains part of a branch of $z.$ If $V$ contains a point of a
branch of $z$ then it contains the entire branch by Lemma
\ref{branch_lem} and hence all branches of $z$.  Thus if $D$ is any
closed disk in $\A_0$ whose boundary consists of branch segments then
the boundary of $D$ is in $V$ and inessential in $\A_0$.  This
boundary cannot be essential in $V$ as $V$ is either simply connected
or annular and essential in $\A_0.$ Hence the boundary of $D$ is
contractible in $V$.  Since $D \subset \A_0$ we conclude that $D
\subset V$.  Since this is true for any $D$ we have shown that $U(z)
\subset V.$
\end{proof}

\begin{prop}\label{separating_prop}
Suppose $z \in \A_0$ is an essential hyperbolic fixed point of $F$.
Then the union of any stable branch and any unstable branch of $z$
contains an essential simple closed curve and a point in each
component of the complement of this curve.
\end{prop}

\begin{proof}
{From} the definition of essential hyperbolic fixed point we know that
there is an essential curve $\beta$ in $\A_0$ consisting of segments
of branches of $z$ which intersect transversely.  The fact that any
stable branch of the fixed point has non-empty transverse
intersection with any unstable branch allows us to deduce that this
curve can be chosen to lie in the union of any pair of branches, one
stable and one unstable.  Let $\alpha_0$ be the frontier of one of the
components of the complement of $\beta$ which intersects one component
of $S^2 \setminus \A_0$ and let $\alpha$ be the frontier of the
component of the complement of $\alpha_0$ which contains the other
component of $S^2 \setminus
\A_0.$ Then $\alpha$ consists of segments of
branches of $z$ which intersect transversely and is simple.  Since it
separates the two components of $S^2 \setminus \A_0$ it is essential.
Since $\alpha$ is in the interior of $U(z)$ the branches of $z$ contain
points from each component of the complement of $\alpha.$
\end{proof}

\bigskip
\section{Regions of instability.}

We begin with the following result which improves Proposition
\ref{per_pt_prop} in the generic case.

\begin{prop}\label{8.1} Suppose $F$ is an area preserving diffeomorphism
of $S^2$ with an invariant annular domain $\A_0$ and $P$ is the set of
periodic points of $F$ with period $q$.  Assume
\begin{enumerate}
\item[i)] 
the set $P$ is empty or finite and all points of $P$ are non-degenerate,

\item[ii)] every elliptic periodic point in $P$ is Moser stable, and
the invariant, $F$ minimal, Jordan curves bounding neighborhoods of
elliptic periodic points have rotation numbers which are not
constant in any neighborhood of the elliptic periodic point, and

\item[iii)] the intersection of any branches of hyperbolic periodic points
in $P$ are transverse and any two branches of the same periodic point have
non-empty intersection.
\end{enumerate}

\noindent
Let $K \subset \A_0$ be an invariant essential continuum. If $p$ is
relatively prime to $q$ and $\rho = p/q$ is in $\rho(K)$, then $P$ is
not empty and there exists an essential hyperbolic periodic orbit in
$K \cap P$ with rotation number $\rho$.  The analogous result holds if
$A = \T^1 \times [0,1]$ and $F:A \to A$ has no periodic points on 
the boundary of $A$ and satisfies i) -- iii).
\end{prop}

\begin{proof}
Taking $F^q$ instead of $F$ and $f^{q} \circ T^{-p}$ instead of $f$ we
can suppose that $\rho=0$.  Replacing $F: \A_0 \to \A_0$ by a lift to
a finite cover we may assume that there are no fixed points except
those which project to fixed points with rotation number $0$.  Note
that a continuum is essential if and only if a lift to a finite cover is
essential.  Hence it suffices to prove the result for this lift since
it has an essential continuum which projects to $K$.  By abuse of
notation we refer to the lift as $F: \A_0 \to \A_0$.

Consider $\widehat F$ the extension of $F$ to $\widehat \A$ the end
compactification of $\A_0$, i.e. the compactification obtained by
collapsing each of the two components of $S^2 \setminus \A_0$ to a
point, obtaining a space $\widehat \A$ homeomorphic to $S^2$.  We
will denote the two points added in this compactification by $z^+$ and
$z^-$.  If $U^+$ and $U^-$ are the components of $\widehat \A
\setminus K$ containing $z^+$ and $z^-$ respectively then they are
simply connected $\widehat F$-invariant domains.  We can apply Theorem
\ref{mather_thm} to show these domains have irrational rotation number
despite the fact that at the points $z^+$ and $z^-$ the homeomorphism
$\widehat F$ is not smooth and neither is the invariant measure.  The
proof of part {\it i)} of \ref{mather_thm} is still valid.

The theory of prime ends allows us to conclude that there is a
closed disk $D^+\subset U^+$ which contains all the fixed points of
$\widehat F$ which are in $U^+$ and such that $i(\widehat
F,D^+)=1$. We can find $D^-\subset U^-$ with the analogous property.

By Proposition \ref{inv_domain_prop} to every inessential hyperbolic
fixed point $z \in \A_0$, we can associate an invariant simply
connected domain $U(z)$.  The same proposition implies there are a
finite set of such domains which are disjoint and which contain all
inessential fixed points of $F$ (all of which, by the remarks above,
can be assumed to have rotation number zero).  Moreover they are
either disjoint from the annular regions $U^+ \setminus \{z^+\}$ and
$U^- \setminus \{z^-\}$ or contained in them.  In each such domain
$U(z)$ which is disjoint from $U^+ \cup U^-$ we may choose a closed
disk $D\subset U(z)$ which contains all the fixed points in $U(z)$ and
such that $i(F,D)=1$. We write $\{D_i\}_{1\leq i\leq k}$, $k\geq 0$,
for the family of disks obtained in this way.

We next consider the finite family $\{z_j\}_{1\leq j\leq
l}$, $l\geq 0$, of fixed points of $F$ which are outside all
the disks $\{D_i\}.$  By the Lefschetz formula, we have
\[
\sum_{j=1}^l i( F,z_j)+\sum_{i=1}^k i(F, D_i) 
+i(\widehat F, D^+)+i(\widehat F,D^-)=2,
\]
and we can write 
\[
\sum_{j=1}^l i(F,z_j)+k=0.
\]

Proposition \ref{per_pt_prop} tells us that one of the integers
$k$ or $l$ is $\geq 1$. We know that $i(F,z_j)=\pm 1$ for every $j$
and that $i(F,z_j)=1$ if $z_j$ is elliptic. We deduce that $l\geq 1$
and that there exists a fixed point $z_j$ such that $i(F,z_j)\leq
0$. This point has to be hyperbolic. It is essential because it is not
contained in the disks $D_i$.

The analogous result for $A = \T^1 \times [0,1]$ follows from similar
but easier reasoning.
\end{proof}

\begin{cor}
Suppose $F$ is an area preserving diffeomorphism of $S^2$ with an
invariant annular domain $\A_0$ and satisfies the hypothesis of
Propositon \ref{8.1} for every $q$.  Then every invariant essential
continuum $K$ in $\A_0$ which is contained in the frontier of a
connected open set $U$ is irrational.
\end{cor}

\begin{proof}
If $K$ is not irrational (see Definition \ref{irr_cont_defn}) it
contains a periodic point $z$.  By Proposition \ref{8.1} we can
suppose that it is an essential hyperbolic periodic point. The set $K$
contains an essential simple closed curve $C$ (in $W^s(z)\cup
W^u(z)$), and a point in each component of the complement of $C$ by
Proposition \ref{separating_prop}. We deduce that the open set $U$
meets the two components of $\A \setminus C$ and doesn't meet $C$. The
contradiction comes from the connectedness of $U$.
\end{proof}

The results above motivate the following definition.

\begin{defn}\label{moserg_def}
Suppose $F$ is an area preserving diffeomorphism
of an annular domain $\A_0$ and $P(q)$ is the set of periodic points of
$F$ with  period $q$.  Assume that for all $q$
\begin{enumerate}
\item[i)] 
the set $P(q)$ is empty or finite and all points of $P$ are non-degenerate,

\item[ii)] every elliptic periodic point in $P(q)$ is Moser stable, and
the invariant, $F$ minimal, Jordan curves bounding neighborhoods of an
elliptic periodic point $z$ have rotation numbers which are not
constant in any neighborhood of $z,$ and

\item[iii)] the intersection of any branches of hyperbolic periodic points
in $P(q)$ are transverse and any two branches of the same periodic point have
non-empty intersection.
\end{enumerate}

Then $F$ will be called {\em Moser generic.} 
\end{defn}

We remark that for $r\geq 3$ the set of Moser generic diffeomorphisms
form a residual subset of ${\rm Diff}^r_{\omega}(\Sigma).$

\begin{prop}
Let $F$ be a Moser generic area preserving diffeomorphism of $S^2$
with an invariant annular domain and let 
$\check F : \check \A \to \check \A$ be the exension to the prime end
compactification of this domain. Then we have the following:

\begin{enumerate}
\item[i)] the union $\mathcal K$ of the irrational invariant
essential continua for $\check F$ is closed;

\item[ii)]the set of the irrational invariant essential
continua for $\check F$ is closed in the Hausdorff topology and the function
$K\mapsto \rho(K)$ is continuous on that set;

\item[iii)]the connected components of $\mathcal K$ are the
maximal irrational invariant essential continua.
\end{enumerate}
\end{prop}

\begin{proof}We first prove assertion {\it ii)}.
Suppose that the sequence $(K_n)_{n\geq 0}$ of irrational invariant
essential continua converges to $K$. This set is an essential
invariant continuum. Let's suppose that it is not irrational. It
contains an essential hyperbolic periodic point $z$, an essential
simple closed curve $C\subset W^s(z)\cup W^u(z)$, a point above that
curve and a point below. For $n$ big enough, $K_n$ also contains a
point in each component of the complement of $C$. We deduce that $K_n$
meets the curve and hence must contain the point $z$. We have a
contradiction to the assumption that $K_n$ is irrational.

Let's suppose that the sequence $(\alpha_n)_{n\geq 0}$ converges to
$\alpha$, where $\alpha_n$ is the rotation number of $K_n$. We can
choose for any $n\geq 0$ an invariant measure $\mu_n\in M(K_n)$. Any
limit point $\mu$ of the sequence $(\mu_n)_{n\geq 0}$ belongs to
$M(K)$ and its rotation number is $\rho(K)$. It is also equal to
$\alpha$. We deduce easily the continuity of the map $K\mapsto
\rho(K)$ from this argument.

Assertion {\it i)} is straightforward. If $(z_n)_{n\geq 0}$ is a sequence in
$\mathcal K$, we can choose for every $n\geq 0$ an irrational
invariant essential continuum $K_n$ which contains $z_n$. Any limit
point of the sequence $(K_n)_{n\geq 0}$ is an irrational invariant
essential continuum which contains $z$.

Assertion {\it iii)} is also very simple. Any connected component
of $\mathcal K$ is closed, invariant, essential and does not contain
any periodic orbit. So it is irrational. Of course, it is maximal
among the irrational invariant essential continua.
\end{proof}

The complement of $\mathcal K$ is an open set $\mathcal U$ which
contains all the periodic orbits of $F$.  Every connected component of
$\mathcal U$ is an invariant essential open annulus.

\begin{defn}\label{region_d}
Suppose $F$ is a Moser generic area preserving diffeomorphism of $S^2$
with an invariant annular domain $\A_0$. Let $\mathcal K$ denote the
union of the irrational invariant essential continua in $\A_0.$
Any invariant component of the complement of  $\mathcal K$ in $\A_0$
will be called a {\it region of instability}.
\end{defn}

Let $U$ be a region of instability and $\check U=U\sqcup C^+\sqcup
C^-$ the prime end compactification of $U$, where $C^+$ is the upper
boundary circle and $C^-$ the lower one.  Let's write $\check F$ for
the extension of $F_{\vert U}$ to $\check U$. We know that the
rotation numbers defined on $C^+$ and on $C^-$ are irrational. Likewise, in
$S^2$ we will denote by $\partial^+U$ and $\partial^-U$ the frontiers
of the two components of the complement of $U$.

\begin{thm}\label{8.6}
Let $F$ be a Moser generic area preserving diffeomorphism of $S^2$
with an invariant annular domain $\A_0$ and let $U$ be a region of
instability in $\A_0.$ We have the following results:

\begin{enumerate}
\item[i)] The rotation set of $\check U$ is a
segment $\rho(U)$ of length $>0$.

\item[ii)] For every rational number $\rho\in\rho (U)$, there
exists an essential hyperbolic periodic orbit of type $(p,q)$ in $U$, where
$\rho=p/q$  and $p$ and $q$ are relatively prime.

\item[iii)] Let $X = K(z)$ denote the closure in $\overline U$ of a
branch of an essential hyperbolic periodic point $z$. Then $X$ is
independent of the essential hyperbolic periodic point $z$ and of the
branch.  It is a subset of any invariant essential continuum contained in
$\overline U$ which is not a subset of $\partial^+U$ or $\partial^-U$.

\item[iv)] The set
$\check X=K(z)$ defined to be the closure in $\check U$ of a branch of $z$ 
is independent of the essential hyperbolic periodic point $z$ and of the
branch. It is a subset of any 
invariant essential continuum contained in $\check U$ which is distinct from
$C^+$ and from $C^-$.
\end{enumerate}
\end{thm}

\begin{proof}First we prove assertions {\it i)} and {\it ii)}.
The set $\overline U$ is an invariant essential continuum which is not
included in $\mathcal K$. So it contains a periodic point. This
periodic point belongs to $U$ because $\partial^+U$ and $\partial^-U$
are irrational invariant essential continua. The rotation numbers of
$C^-$ and $C^+$ being irrational, the set $\rho(U)$ is a segment of
length $>0$. The extension of $F$ to the prime-end compactification of
$U$ satisfies the hypothesis given at the beginning of Paragraph 7 and
we can apply Proposition \ref{8.1} to the essential continuum $\check
U$ itself: for every rational number $\rho\in\rho(U)$ there exists in
$\check U$ an essential hyperbolic periodic point of type $(p,q)$,
where $\rho=p/q$ and $p$ and $q$ are relatively prime.  But this point
belongs to $U$ and not to the boundary curves.

We now prove assertion {\it iii)}. Let $K$ be an essential
invariant continuum contained in $\overline U$ and not a subset of 
$\partial^+U \cup \partial^-U$. The boundary of $U_{K}^+$ is an
irrational invariant essential continuum. It cannot meet $U$ (since it
is irrational) and must be contained in $\partial U$.  More precisely
it is contained in $\partial^+U$ because it is connected and because
$K$ is not a subset of $\partial^-U$.  In the same way the boundary of
$U_{K}^-$ is contained in $\partial^-U$. Hence we conclude that $K$
meets the two boundary components of $U$. If $z$ is an essential
hyperbolic periodic point, we can find an essential curve contained in
$W^s(z)\cup W^u(z)$. The set $K$ meets that curve and must contain
$K(z)$. In particular, if $z'$ is another essential hyperbolic
periodic point contained in $U$ the set $K(z')$ contains $K(z)$.

The proof of the assertion {\it iv)} is exactly the same.
\end{proof}

We deduce immediately the following result about periodic points:

\begin{cor} \label{8.7}
Let $z$ and $z'$ be two essential hyperbolic periodic points of
$F$. The following conditions are equivalent:

\begin{enumerate}
\item[i)] the points $z$ and $z'$ are in the same
region of instability;

\item[ii)] we have $K(z)=K(z')$;

\item[iii)] we have $W^u(z)\cap W^s (z')\not=\emptyset$.
\end{enumerate}
\end{cor}

We next consider a result that generalies a property of 
regions of instability of twist maps:

\begin{thm}\label{8.8}
Let $U$ be a region of instability. For any neighborhood $W$ of
$\partial^+U $ in $\overline U$ whose closure doesn't contain $\partial^-U$,
the connected component
of $\displaystyle\bigcap_{k\in\Z}F^{-k}(\overline W)$ which contains
$\partial^+ U$ is equal to $\partial^+ U$.
\end{thm}

\begin{proof}
The connected component $K$ of
$\displaystyle\bigcap_{k\in\Z}F^{-k}(\overline W)$ which contains
$\partial^+ U$ is an invariant essential continuum and the boundary of
$U_{K}^-$ is an irrational invariant essential continuum. This
boundary is connected, doesn't meet $U$ and doesn't contain
$\partial^- U$: it is contained in $\partial^+U$. We deduce that
$K=\partial^+U$.
\end{proof}

We obtain the following corollaries:

\begin{cor} \label{8.9}
Let $U$ be a region of instability.  In any neighborhood $W$ of
$\partial^+U$ in $\overline U$ there exists a periodic orbit.
Also in $W$ there is a homoclinic point of any essential periodic
point in $U$.
\end{cor}

\begin{proof}
We write $z^+$ for the upper end of $U$ and $z^-$ for the lower
end. Let's consider the extension $\widehat F$ of $F_{\vert U}$ to the
end compactification $\widehat U$. The end $z^+$ is a fixed point of
$\widehat F$ and the sequence $\left(i(\widehat F^k,z^+)\right)_{k\geq 1}$ of
Lefschetz indices is well defined and constant equal to $1$. Indeed,
the rotation number of $C^+$ is irrational. Let's consider a closed
disk $D$ which contains $z^+$ in its interior and which is contained
in $(W \setminus \partial U)\sqcup\{z^+\}$. The connected component
of $X=\displaystyle\bigcap_{k\in\Z}\widehat F^k(D)$ which contains $z^+$
is equal to $z^+$ by Theorem \ref{8.8}. So, we can find a closed disk
$D'\subset D$ which contains $z^+$ in its interior and such that
$\partial D'\cap X=\emptyset$. We deduce that the set
$\displaystyle\bigcap_{k\in\Z}\widehat F^k(D')$, which is contained in
$X$, doesn't meet the boundary of $D'$. Such a set is called an
{isolating block}. To such a set we can associate a sequence of
integers
\[
a_k=\left(i\left(\widehat F^k,\bigcap_{0\leq i\leq
k}f^{-i}(D')\right)\right)_{k\geq 1}
\]
 and this sequence is non
positive for infinitely many values of $k$ (see \cite {F2} or \cite {LY}). The
integer $a_k$ is equal to the sum of the Lefschetz indices (for
$\widehat F^k$) of all the fixed points of $\widehat F^k$ whose $\widehat F$ 
orbit lies entirely in $D'$ if there are a finite number of such points. For
any value $k$ such that $a_k\leq 0$ there exists a fixed point of
$\widehat F^k$ distinct from $z^+$, whose orbit lies entirely in $D'$.

To prove that there is a homoclinic point in $W$ we consider an
arbitrary essential fixed point $z$ of $F^q$ and let $C$ be an
essential simple closed curve in $U$ consisting of segments of
$W^s(z)$ and $W^u(z)$.  We define
\[
B_n = \bigcup_{k=-n}^n F^{-k}(C)
\]
and let $A_n$ be the closure of the component of the complement of
$B_n$ in $\overline U$ which contains $\partial^+U.$ We note that
$F(B_n) \subset B_{n+1}$ so $A_{n+1} \subset F(A_n).$ {From} this it
follows that $\cap_{n \ge 0} A_n \subset \cap_{n \ge 0}F( A_n) =
F(\cap_{n \ge 0} A_n)$ and hence that $X = \cap_{n \ge 0} A_n$ is
forward invariant under $F^{-1}$ (and similarly under $F).$ We
conclude $X$ is invariant and from Theorem \ref{8.8} that $X =
\partial^+U$ and since the family $\{A_n\}$ is a nested sequence of
compact sets we know that $A_n \subset W$ for $n$ sufficiently large.
The frontier of $A_n$ consists of finitely many segments of $W^s(z)$
and $W^u(z)$ and hence must contain a homoclinic point of $z.$
\end{proof}

\begin{rem}
In the case where $F$ is a twist map, the Birkhoff theory tells us
that any irrational invariant essential continuum is the graph of
Lipschitz map $\psi: \T^1\to[0,1]$. In particular we don't need the
prime end theory. We have a natural order on the set of these curves
($K<K'$ if $K'$ is above $K$) and the function $K\mapsto\rho(K)$ is
strictly monotone. There are some other results that we know in the
case of a twist map.
\end{rem}

\medskip
Using variational methods, Mather \cite {M2} proved that there exists
a point $z\in U$ whose $\alpha$-limit set is contained in
$\partial^+U$ and whose $\omega$-limit set is contained in
$\partial^-U$, if $U$ is a region of instability of a twist map.
There exists a topological proof of this result
\cite {L} which uses the Birkhoff theory for positively and negatively
invariant essential continua.   We now prove that this result remains
valid for generic area preserving diffeomorphisms without the requirement
of the twist condition.

\begin{prop}\label{connecting_t} 
Let $U$ be a region of instability.
There exists a point $x$ in $U$ such that $\omega(x) \subset 
\partial^+U$ and $\alpha(x) \subset \partial^-U$.
\end{prop}

\begin{proof}

If $C$ is an essential simple closed curve in $U$ and if $N$ is the
component of its complement (in $\overline U$) which contains $\partial^+
U$, a well-known result of Birkhoff (see \cite{Bi2} Chapter VIII, section 8)
implies that the component of the compact set
\[
\bigcap_{k=0}^{\infty} F^{-k}(N)
\]
containing $\partial^+U$, must meet $C$. This component, $X$, is a
positively invariant set containing the set $\partial^+U$ but not
equal to it.

Theorem \ref{8.8} tells us that the globally invariant set
$\cap_{n\geq 0} F^n(X)$ is reduced to $\partial^+U$. So for any point
$w\in X$ we have $\omega(w)\subset\partial ^+U$. Theorem \ref{8.8}
tells us that the connected globally invariant set 
$cl(\cup_{n\geq 0} F^{-n}(X))$ meets $\partial^-U$.  If $z\in U$ is
an essential hyperbolic periodic point, Proposition
\ref{separating_prop} says there exists an essential closed curve
contained in the union of any unstable branch of $z$ and any stable
branch. If $n$ is large enough, the set $F^{-n}(X)$ must meet this
curve. We deduce that $X$ meets both unstable branches of $z$. More
precisely, there exists an open interval $\gamma$ of $W^u(z),$ which
contains $z,$ is disjoint from $X,$ and whose boundary points belong
to $X$. This arc is contained in a component $V$ of $S^2\setminus X$
and separates $V$ into two components (because $X$ is connected) which
are connected components of $S^2\setminus (\gamma\cup X)$.

Similarly we can construct a negatively invariant compact set $Y$ in
$U\cup\partial^- U$, containing the set $\partial^-U$ but not equal to it.
For any point $w\in Y$ we have $\alpha(y)\subset\partial^-U$. 
Also as before, $Y$ meets both branches of $W^s(z)$.

Consider two points $w_1,w_2\in Y$, one in each branch of $W^s(z).$ 
For large $n$ the points $F^n(w_1),F^n(w_2)Y$ will be 
close to $z$, one in each component of $V\setminus \gamma$.
Hence the connected set $F^n(Y)$ meets the two components of $V\setminus
\gamma$. The set $F^n(Y)$ doesn't intersect $\gamma$, so it intersects $X$.
Any point of intersection of $X$ and $F^n(Y)$ satisfies the proposition.
\end{proof}

Among the periodic orbits that we can find in a region of instability
there are particular orbits, called {\it well-ordered orbits} or {\it
Birkhoff orbits}. Assertion {\it ii)} of Theorem \ref{8.6} can be
improved: the periodic orbit can be chosen among the Birkhoff
orbits. There is a natural generalization of Birkhoff orbits in the
general case which is the notion of {\it unlinked orbit}. Boyland
\cite {Bo} proved that for a homeomorphism of the closed annulus the
existence of a periodic orbit of type $(p,q)$ implies the existence of
an unlinked periodic orbit of the same type.  Moreover there is a
Lefschetz formula for these orbits. With these remarks it is not
difficult to improve assertion {\it ii)} of Theorem \ref{8.6}: the
periodic orbit can be chosen among the unlinked orbits.

In the case of a twist map there exists in a region of instability $U$
a well-ordered invariant Cantor set of rotation number $\rho$ for
every irrational number in the interior of the rotation set of
$U$. Such a set is called an {\it Aubry-Mather} set. Using symbolic
dynamics and the heteroclinic intersections between periodic points of
rotation number smaller and bigger than $\rho$, is it possible to find
an invariant (or an unlinked invariant) compact set of rotation number
$\rho$ in the general case (see \cite {H}) ?

The following result is closely related to a result of
Boyland \cite{Bo}.

\begin{prop}
Let $F$ be a Moser generic area preserving diffeomorphism of $S^2$
with an invariant annular domain $\A_0$ and let $K$ be an essential
continuum which is the closure of one branch of a hyperbolic fixed
point $z$ of $F$.  
Let $\widetilde z$ be a point in $\pi^{-1}(z)$ in the universal cover
$\widetilde \A$ of $\A_0$ and let $f:\widetilde \A \to \widetilde \A$
be the lift which fixes $\widetilde z$.
If $\widetilde K \subset \widetilde A$ is
$\pi^{-1}(K)$ and $\rho( \widetilde z, f) \in int \rho( \widetilde K, f)$
then $\widetilde K$ is the closure of any branch of
$\widetilde z$, i.e. each branch of $\widetilde z$ in $\widetilde A$
is dense in $\widetilde K$.
\end{prop}

\begin{proof}
Let $\Gamma$ be a branch of $\widetilde z$ for the lift $f.$ We can
assume without loss of generality that it is an unstable branch.
Since $z$ is essential there is a $k \ne 0$ such that $\Gamma$
intersects $W^s(T^k(\widetilde z))$ transversely.  The $\lambda$-lemma
then implies the closure of $\Gamma$ contains the closure of
$W^u(T^k(\widetilde z))$.  The $\lambda$-lemma also implies that the
set $J$ consisting of all values of $k$ for which $\Gamma$ intersects
$W^s(T^k(\widetilde z))$ is closed under addition.

Another straight-forward $\lambda$-lemma argument and the fact that
the rotation number of $z$ is in the interior of the rotation set
for $K$ shows that $J$ contains both positive and negative elements.
{From} this we may conclude that $J$ is a subgroup of $\Z,$ say the 
subgroup $d\Z.$ 

If we replace $\A_0$ and $F$ by a $d$-fold cover of $\A_0$ and a lift
of $F$ to it then clearly nothing has changed in what we want to prove
except that now $d$ has been replaced by $1$.  I.e. we may assume
without loss of generality that $d = 1$ or that $J = \Z.$
Hence we may assume that the closure of $\Gamma$ contains the closure of
$W^u(T^k(\widetilde z))$ for every $k$.

Now let $\widetilde y$ be an element of $\widetilde K$ and let $y =
\pi(\widetilde y).$ The fact that $W^u(z,F)$ is dense in $K$ implies
that for any $\ep > 0$ there is a point $x$ of $W^u(z,F)$ which
is within $\ep$ of $y$.  Equivalently, there is an integer $k$ and
a point $\widetilde x$ of $W^u(T^k(\widetilde z,f))$ which is within
$\ep$ of $\widetilde y.$ Hence $\Gamma$ is dense in $\widetilde K.$
\end{proof}

\section{Generic diffeomorphisms of the sphere}

In this section we consider a $C^r$-diffeomorphism $F$ of $S^2$ which
preserves a volume form $\omega$.
We will further assume $F$ is Moser generic (see \ref{moserg_def}).

For every integer $q\geq1$ we write ${\rm Hyp}_q(F)$ (resp. ${\rm
Ell}_q(F)$) for the set of hyperbolic (resp. elliptic) fixed point of
$F^q$. We note
\[
{\rm Hyp}(F)=\bigcup_{q\geq 1}{\rm Hyp}_q(F)
\]
the set of hyperbolic (resp. elliptic) periodic points and
\[
{\rm Ell}(F)=\bigcup_{q\geq 1}{\rm Ell}_q(F)
\]
the set of elliptic periodic points.

For every integer $q\geq 1$, we note $\K_q$ the set of non empty
continua which are (globally) invariant by $F^q$ and
$\K_q^*\subset\K_q$ the set of continua which do not consist of a single
fixed point of $F^q$. We note $\K=\cup_{q\geq 1} \K_q$ the sets of
periodic continua and $\K^*=\cup_{q\geq 1} \K_q^*$ the set of non
trivial periodic continua.

We recall a property we will use in this section -- that a closed set
$K\subset S^2$ is connected if and only if the complement $U=S^2
\setminus K$ is simply connected.

{From} our earlier results we have the two following properties

\begin{prop} \label{9.1}
Any continuum $K\in\K^*$ which contains a hyperbolic periodic point
$z$ contains the four branches of $z$.
\end{prop}

\begin{prop} \label{9.2}
Any continuum $K\in\K_q$ which doesn't contain a periodic point
separates the sphere in two connected and simply connected components which are
$F^q$-invariant. Moreover the rotation set of the continuum consists
of a single irrational number (defined up to sign).
\end{prop}

A periodic continuum without any periodic orbit will be called a {\it
periodic irrational continuum} and we will define its {\it rotation
number} which is an element of $\T^1$ defined up to sign.

Let $K$ and $K'$ be two non disjoint periodic irrational continua and
write $U$ and $V$ for the components of $S^2 \setminus K$ and $U'$ and
$V'$ for the components of $S^2 \setminus K'$.  The set $K\cup K'$ is
a periodic irrational continuum and the set 
\[
S^2 \setminus (K\cup K')=(U\cap U')\sqcup(U\cap V')\sqcup(V\cap
U')\sqcup(V'\cap U')
\]
has exactly two connected components by \ref{9.2}.  Each set $U$ and
$V$ contains a periodic point and cannot be contained in $K'$: we
deduce that one of the sets $U\cap U'$ or $U\cap V'$ is non empty as
is one of the sets $V\cap U'$ and $V\cap V'$. Doing the same thing for
$U'$ and $V'$ we can suppose (exchanging $U$ and $V$ if necessary)
that the two components of $S^2 \setminus K\cup K'$ are $U\cap U'$ and
$V\cap V'$ and that $U\cap V'=V\cap U'=\emptyset$. The sets $U$ and
$U'$ contain the same periodic points as do the sets $V$ and $V'$. The
open sets $U\cup U'$ and $V\cup V'$ are connected and simply
connected, because $U\cap U'$ and $V\cap V'$ are non-empty and
connected.  Thus the sets $U\cup U'$ and $V\cup V'$ are the components
of the complement of the irrational continuum $K\cap K'$.

By induction we can prove the analogous result for more than two
periodic irrational continua.  More precisely, if $(K_i)_{1\leq i\leq
n}$ is a family of periodic irrational continua such that the union is
connected, the union is a periodic irrational continuum.  It's
complement has two components $\cap_{1\leq i\leq n}U_i$ and
$\cap_{1\leq i\leq n} V_i$ where $U_i$ and $V_i$ are the components
(well chosen) of $S^2 \setminus K_i$. The set $\cap_{1\leq i\leq n}
K_i$ is a periodic irrational continua and the components of its
complement are $\cup_{1\leq i\leq n}U_i$ and $\cup_{1\leq i\leq n}
V_i$.

Let $K$ be a periodic irrational continuum whose interior $O$ is not
empty and consider the two components $U$ and $V$ of $S^2 \setminus
K$. The periodic irrational continua $\partial U$ and $\partial V$
don't intersect.  This is because the sets $U$, $V$ and $O$ are in
three separate components of $S^2 \setminus (\partial U\cup\partial
V)$ and if $\partial U$ meets $\partial V$ their union is an invariant
irrational continuum whose complement has only two components.  Thus
the set $O$ must be an open annulus.

Let $\Gamma^s$ be a stable branch of $z\in {\rm Hyp}_q(F)$ and $\Gamma
^u$ an unstable one. Then, the two branches intersect, they have the
same closure $K(z)\in\K_q$ and there is a fundamental system of
neighborhoods of $z$ which are Jordan domains with a boundary
contained in $\Gamma^s\cup\Gamma^u$.  We know that the equality
$K(z)=K(z')$ occurs if and only if any stable branch of $z$ meets any
unstable branch of $z'$.  We know that there is no strict inclusion
$K(z)\subset K(z')$.

Let's prove the following fundamental result.

\begin{prop}\label{union_prop}
If $(U_i)_{i\in I}$ is a family of periodic simply connected open
sets, the union $U=\cup_{i\in I} U_i$ is simply connected.
\end{prop}

\begin{proof}
It is sufficient to prove this result in the case where $I$ is finite
and the $U_i$ are connected. In fact, appealing to induction, it is
sufficient to prove that the union of two invariant simply connected
domains $U$ and $U'$ is simply connected.

The proof is obvious if one of the complements $S^2-U$ or $S^2-U'$ is
empty or reduced to a point. We will suppose that none of these
properties occurs. The boundaries $\partial U$ and $\partial U'$ are
irrational invariant continua and we can write
\[
S^2=U\sqcup \partial U\sqcup V=U'\sqcup \partial U'\sqcup V'
\]
 where $V$ and $V'$ are some
invariant simply connected domains.

If the sets $\partial U$ and $\partial U'$ meet, their union is an
irrational invariant continuum and the complement 
\[
S^2 \setminus (\partial U \cup \partial U')=(U\cap U')\sqcup(U\cap
V')\sqcup(V\cap U')\sqcup(V\cap V')
\]
 has exactly two components, by \ref{9.2}.  At
least one of the sets $V\cap U'$ or $V\cap V'$ is not empty (because
$\partial U'$ has no interior and cannot contain $V$). We deduce that
$U\cap U'$ is either empty or connected and that $U\cup U'$ is simply
connected.

If $\partial U$ is contained in $U'$, then one of the open sets $U$ or
$V$ doesn't meet the connected set $S^2 \setminus U'$ and is contained
in $U'$. In the first case $U\cup U'=U'$ is simply connected, in the
second case we have
\[
S^2 \setminus U=\partial U\cup V\subset U'
\]
so $U\cup U'= U$ is simply connected.

If $\partial U$ is contained in $V'$, then one of the open sets $U$ or
$V$ doesn't meet the connected set $S^2 \setminus V'$ and is contained
in $V'$. In the first case $U$ and $U'$ are disjoint and $U\cup U'$ is
simply connected, in the second case we have 
\[
S^2 \setminus U=\partial U\cup V\subset V'\subset S^2 \setminus U'
\]
and $U\cup U'=S^2$ is simply connected.
\end{proof}

\begin{cor} \label{9.4}
Let $U$ and $U'$ be two disjoint periodic simply connected domains. If
$(U_i)_{i\in I}$ is the family of all periodic domains disjoint from $U$
and $U'$, then the complement $K(U,U')$ of 
\[
U\cup U'\cup_{i\in I} U_i
\]
is the smallest periodic continuum such that $U$ and $U'$ are
components of its complement. Moreover $K(U,U')$ is invariant if 
$U$ and $U'$ are invariant.
\end{cor}

\begin{proof}
By Proposition \ref{union_prop}, we know that $K(U,U')$ is
connected. It is clear that $U$ and $U'$ are connected components of
$S^2 \setminus K(U,U')$. If $K$ is a periodic continuum such that $U$
and $U'$ are components of $S^2 \setminus K$, then the other
components are in the family $(U_i)_{i\in I}$, so $K$ contains
$K(U,U')$. To see that $K(U, U')$ is invariant if $U$ and $U'$ are we
observe that $\cup_{i\in I}U_i$ is $F$-invariant in that case.  We
still have to prove that $K$ is periodic, but this follows by
considering $F^q$ where $q$ is a common period for $U$ and $U'$ and
noting that $K(U,U')$ is the same for $F$ and $F^q$, so our previous
remark implies $K(U,U') = F^q(K(U,U')).$
\end{proof}

\begin{prop} \label{9.5_prop}
Let $U$ and $U'$ be two disjoint periodic simply connected domains. If
$S^2 \setminus (U\cup U')$ contains a periodic point, so does
$K(U,U')$. Moreover every periodic point $z\in K(U,U')$ is hyperbolic.
\end{prop}

\begin{proof}
We can suppose that $U$ and $U'$ are invariant. If $K(U,U')$ has no
periodic point it is an irrational continuum, so the components of
$K(U,U')$ are exactly $U$ and $U'$ and we deduce that $K(U,U')=S^2
\setminus (U\cup U')$.

Let's prove now the second part of the proposition. There is no
periodic point in $\partial U$ and in $\partial {U'}$ because these
two sets are periodic irrational continua. So every elliptic periodic
point $z\in S^2 \setminus (U\cup U')$ has a periodic simply connected
neighborhood, bounded by a Jordan curve, which meets neither $U$ nor
$U'$. We deduce that $z\not\in K(U,U')$.
\end{proof}

We can improve this result.

\begin{prop}
Let $U$ and $U'$ be two disjoint invariant simply connected
domains. Consider a fixed point $z\in U$, a fixed point $z'\in U'$ and
a lift $f$ of $F_{\vert S^2 \setminus \{z,z'\}}$ to the universal
covering space of the annulus $S^2 \setminus \{z,z'\}$. Let $p\in\Z$
and $q\geq 1$ relatively prime. If $S^2 \setminus (U\cup U')$ contains
a periodic point of type $(p,q)$, so does $K(U,U')$.
\end{prop}

\begin{proof}
The set $Z$ of periodic points of type $(p,q)$ in $S^2 \setminus
(U\cup U')$ is finite. If it is disjoint from $K(U,U')$, there exists
a smallest finite family $(U_i)_{1\leq i\leq n}$ of components of $S^2
\setminus K(U,U')$ which covers $Z$. Every $U_i$ is $F^q$-invariant
and every fixed point of $F^q$ in $U_i$ belongs to $Z$. The
contradiction comes from the Lefschetz formula
\[
0=\sum_{z\in Z} i(F^q,z)=\sum_{i=1}^n i(F^p,U_i)=n.
\]

\end{proof}

We can be more precise.

\begin{prop} \label{9.7_prop}
Let $U$ and $U'$ two disjoint invariant simply connected domains. For
any periodic hyperbolic point $z\in K(U,U')$, the set $K(z)$ (defined to
be the closure of a branch of $z$) is
invariant and separates $U$ and $U'$.  Moreover if $\Gamma^s$ is a
stable branch of $z$ and $\Gamma ^u$ and unstable one, there is a
Jordan curve which separates $U$ and $U'$ and which is contained in
$\Gamma^u\cup \Gamma^s$.
\end{prop}

\begin{proof}
The periodic continuum $K(z)$ is contained in $K(U,U')$. If it doesn't
separate $U$ and $U'$, the component of $S^2 \setminus K(z)$ which
contains $U$ and $U'$ is a periodic simply connected domain $V$ and
$z\not\in\overline V$ since $z$ has a fundamental system of
neighborhoods whose boundaries are in the branches of $z$.  We deduce
that the component of $S^2 \setminus \overline V$ which contains $z$
is a periodic domain disjoint from $U$ and $U'$. It is disjoint from
$K(U,U')$ by definition of this set so we have a contradiction.

Let's prove now that $K(z)$ is invariant by showing that the orbit of
$z$ is contained in $ K(z)$. Consider the components $V$ and $V'$ of
$S^2 \setminus K(z)$ which contain respectively $U$ and $U'$. The set
$\partial V$ is an $F^q$-invariant irrational continuum and meets its
image because it separates $U$ and $U'$ and $F$ is area
preserving. We deduce that $K = \cup_{0\leq k<q}F^k(\partial V)$ is an
$F$-invariant irrational continuum. The orbit of $z$ doesn't meet the
component $\cap_{0\leq k<q}F^k(V)$ of $S^2 \setminus K$ which contains
$U$, it is contained in the other component which is disjoint from
$V$. Similarly the orbit doesn't meet $V'$.  It doesn't meet any
other component of $S^2 \setminus K(z)$ because it is contained in
$K(U,U')$. We have proven that it is contained in $K(z)$.

We deduce that any stable branch of $z'$ in the orbit of $z$ meets any
unstable branch of $z$. We deduce from this fact and from the $\lambda$-lemma 
the end of the proposition in the case where $z$ is a periodic
point of period $q\geq 2$.  But if $z$ is a fixed point then we choose
a hyperbolic periodic point $w \in K(z)$ of period at least two and
observe that $K(z) = K(w).$   The fact that $W^s(w)$ intersects
$\Gamma^u$ transversely and  $W^u(w)$ intersects
$\Gamma^s$ transversely allows us to apply the $\lambda$-lemma 
once again and obtain a simple Jordan curve separating $U$ and $U'$
which is a subset of $\Gamma^u \cup \Gamma^s.$

\end{proof}

\begin{cor} \label{9.8_cor}
Any non-trivial invariant continuum $K$ which is not irrational contains
an invariant continuum $K(z)$, for some $z\in {\rm Hyp}(F)$.
\end{cor}

\begin{proof}
If the complement of $K$ has at least two invariant domains, this is a
consequence of Propositions \ref{9.5_prop} and \ref{9.7_prop}. In the
other case the Lefschetz formula tells us that there exists at least
one fixed point in $K$. If $K$ contains a hyperbolic fixed point, the
result is obvious. If $K$ contains an elliptic fixed point $z$, we can
find two minimal invariant curves $C_1$ and $C_2$ of different
irrational rotation number surrounding $z$ in a neighborhood of that
point. These curves meet $K$ and are contained in it. By
the Poincar\'e-Birkhoff theorem, we can apply Proposition 9.5 and 9.7 to
the domain $U$ containing $z$ bounded by the closest curve and to the
domain $U'$ which doesn't contain $z$ and bounded by the other
curve. There exists $z'\in{\rm Hyp}(F)$, such that $K(z')$ is invariant
and a simple curve separating $U$ and $U'$ and contained in
$W^u(z')\cup W^s (z')$. This curve meets the connected set $K$ and we
have $K(z')\subset K$.
\end{proof}

\begin{cor}
The invariant continua $K(z)$, $z\in {\rm Hyp}(F)$, are the minimal
elements in the set of the non trivial invariant continua which are
not irrational.
\end{cor}

\begin{proof}
This is an immediate consequence of Corollary \ref{9.8_cor} and of the
fact that there is no strict inclusion $K(z)\subset K(z')$, for
$z,z'\in{\rm Hyp}(F)$.
\end{proof}

\begin{prop} \label{9.10_prop}
The set of irrational invariant continua and continua $\{z\}$,
$z\in{\rm Ell}_1(F)$, is closed in the Hausdorff topology. The union
of all these sets is closed.
\end{prop}

\begin{proof}
Recall that the set $\K_1$ consists of invariant continua.  The result
is a consequence of Corollary \ref{9.8_cor} and of the following
facts:

\begin{enumerate}
\item[$\bullet$] the set $\K_1$ is closed in the Hausdorff topology;

\item[$\bullet$] any continuum $\{z\}$, $z\in{\rm Hyp}_1(F)$, is
isolated in $\K_1$;

\item[$\bullet$] for any $z\in {\rm Hyp}(F)$, the set of continua
$K\in\K_1$ which contain $K(z)$ is open in $\K_1$.
\end{enumerate}
\end{proof}

\begin{rem}
By Birkhoff's theory, if $z$ is a elliptic fixed point and $V$ a small
neighborhood of $z$, there is a neighborhood of $\{z\}$ in the set
defined in Proposition \ref{9.10_prop} which consists of $\{z\}$ and of
the invariant curves surrounding $z$ and contained in $V$. If $z'$ is
another fixed point, the set which consists of $\{z\}$, $\{z'\}$ and
of the irrational invariant continua which separate $z$ and $z'$ is
closed.
\end{rem}

\begin{prop}\label{9.12_prop}
The set consisting of irrational invariant continua, of invariant sets
$K(z)$, for $z\in{\rm Hyp}(F)$, and of continua $\{z\}$, $z\in{\rm
Ell}_1(F)$, is closed in the Hausdorff topology.
\end{prop}

\begin{proof}
Indeed, any convergent sequence $K(z_n)$, $z_n\in{\rm Hyp}(F)$, of
invariant continua converges either to a trivial continuum $\{z\}$,
$z\in{\rm Ell}_1(F)$, or to an irrational invariant continuum,
or to an invariant continuum which contains an invariant set
$K(z)$, $z\in {\rm Hyp}(F)$. In the last situation the set $K(z)$ is
contained in $K(z_n)$, for $n\geq n_0$ big enough. The sequence
$K(z_n)_{n\geq n_0}$ is constant and equal to $K(z)$.
\end{proof}

\bigskip
A periodic domain $U$ of period $q\geq 1$ will be called a {\it
strictly periodic} domain if the domains $F^{k}(U)$, $0\leq k<q$ are
all disjoint. Note that the closures also are disjoint. Indeed if
$\partial U$ meets $\partial U'$ where $U' = F^k(U)$, $0<k<q$, the
union of these two sets is an $F^q$-invariant irrational continuum.  Writing
$V=S^2 \setminus \overline U$ and $V' = S^2 \setminus \overline {U'}$,
the components of the complement must be $U\cap V'=U$ and $U'\cap
V=U'$ (because $U$ and $U'$ are disjoint) and we must have
$V \cap V' = \emptyset.$  This is possible only if
$q=2$, and $F(U) = U', F(U') = U$  and $S^2 = \overline U \cup \overline {U'}.$
But that leads to a contradiction since there are no fixed point in 
$U \cup U'$ there must (by the Lefschetz theorem) be a fixed point in
$\partial U \cup \partial U'$.

\begin{prop} \label{9.13_prop}
The union $X_1$ of the irrational invariant continua, of the invariant
sets $K(z)$, $z\in{\rm Hyp}(F)$ and of the elliptic fixed points is an
invariant continuum. Its complement is the union of the simply
connected strictly periodic domains of period $q\geq 2$
\end{prop}

\begin{proof}
The fact that it is closed is a consequence of Proposition
\ref{9.12_prop}.  We prove that it is connected by showing that any
component $U$ of $S^2 \setminus X_1$ is simply connected.

We first observe that $U$ is not invariant. If it were invariant, its
closure would be an invariant continuum which is not irrational
(because $U\not\subset X_1$). It would contain an invariant continuum
$K(z)$, $z\in{\rm Hyp}(F)$, which is disjoint from $U$ (because $U\cap
X_1=\emptyset$). A contradiction arises from the fact that there is a
fundamental system of neighborhoods of $z$ whose boundary is in
$W^s(z)\cup W^u(z)$.

Next we prove that $U$ is simply connected. Any Jordan curve $\gamma$
contained in $U$ is disjoint from its image.  The map $F$ being
volume-preserving, one of the two domains bounded by $\gamma$ is
disjoint from its image. Any component of $X_1$ being invariant, we
deduce that this domain is disjoint from $X_1$ and contained in $U$.
We conclude that $U$ is simply connected.  We have proven that any
component of $S^2 \setminus X_1$ is a strictly periodic simply
connected domain of period $q\geq 2$.

Conversely, a strictly periodic domain $U$ of period $q\geq 2$ doesn't
contain any fixed points. It cannot meet any fixed continuum $K(z)$,
$z\in{\rm Hyp}(F)$ since otherwise it will meet $W^s(z)$ and the
closure $\overline U$ which is disjoint from its image will contain
$K(z)$.

Finally we must show that $U$ cannot meet any irrational invariant
continuum $K$.  If it did then we could not have $K \subset U$ so $K
\cap \partial U \neq \emptyset.$ It follows that $\partial U$ is a
periodic irrational continuum and by our remarks above that $K' = K
\cup \partial U \cup F(\partial U)$ is a periodic irrational
continuum.  Since $K$ has no interior, the complement of $K'$ would
necessarily have at least $3$ components since it would contain a
point from $U$, a point from $F(U)$ and a point in neither of them.
Each of these would be in different components of $S^2\setminus K'.$
This contradicts \ref{9.2}.
\end{proof}

We will write $E$ for the union of the periodic irrational continua and
of the elliptic periodic points, $H$ for the union of the $K(z)$,
$z\in{\rm Hyp}(F)$ and $I$ for the complement of $E\cup H$. So we have
a decomposition 
\[
S^2=E\cup H\sqcup I.
\]

\begin{prop} \label{9.14_prop}
We have the following properties~:

\begin{enumerate}
\item[i)]the sets $E$ and $H$ are $F_{\sigma}$ ;

\item[ii)] we have $I\subset\partial E=\partial H$ ;

\item[iii)] the set $E\cup H$ is connected and dense~;

\item[iv)] a point $z$ belongs to $I$ if and only if there
exists a decreasing sequence $(U_n)_{n\geq 0}$ of strictly periodic
simply connected domains which contain $z$ such that the sequence of
periods tends to the infinity.
\end{enumerate}
\end{prop}

\begin{proof}
We can write $E=\cup_{q\geq 1} E_q$ and $H=\cup_{q\geq 1} H_q$, where
$E_q$ is the union of the $F^q$-invariant irrational continua and of
the elliptic fixed points of $F^q$ and $H_q$ the union of the sets
$K(z)$, $z\in{\rm Hyp}(F)$, which are $F^q$-invariant. Each set $E_q$
is closed and each set $H_q$ is an $F_{\sigma}$, moreover $E_q\cup H_q$
is closed, connected and contains $H_1\cup E_1$.  The first assertion
is clear as is the connectedness of $E\cup H$.  The last assertion is an
immediate consequence of Proposition \ref{9.13_prop}. The density of
$E\cup H$ follows from the inclusion
$I\subset\overline E$ that we will prove now. For any decreasing
sequence $(U_n)_{n\geq 0}$ of strictly periodic domains which contains
$z\in I$ such that the sequence of periods tends to infinity, the
volume of $U_n$ tend to zero. So $\partial U_n$ which is contained in
$E$ meets a given neighborhood of $z$ if $n$ is big
enough.

Since $I\subset\overline E$ we deduce that $\partial H \subset
\partial E.$ If we show $I\subset\overline H$ then it follows that
$\partial E \subset \partial H$ and we can conclude assertion {\it
ii)}.  If $z_0$ is in $I$ then by {\it iv)} $z \in U$ where $U$ is a
strictly periodic simply connected domain of high period (and hence of
small measure).  Given $\epsilon>0$ we may choose $U$ so that any
point of $U$ is within $\epsilon$ of $\partial U.$ We also know that
$z \in U' \subset U$ where $U'$ is a strictly periodic simply
connected domain of higher period than $U.$ The annular domain $U
\setminus \overline U'$ is invariant under some iterate of $F$ and
contains periodic points, so we may apply Proposition \ref{8.1} to its
prime end compactification to conclude that it has an essential
periodic point $z_0$.  Then $K(z_0)$ separates $z$ and $\partial U.$
It follows that there must be a point of $K(z_0) \subset H$ within
$\epsilon$ of $z$ and hence $I\subset\overline H$.
\end{proof}

\begin{thm}\label{dense_stable}
There exists a dense $G_{\delta}$, $\G$ in ${\rm Diff}_{\omega}^r(S^2)$
such that the interior of $E$ is empty if $F\in \G$. The union of
stable manifolds of hyperbolic periodic points is dense in $S^2$
for $F$ in $\G$. 
\end{thm}

\begin{proof}
It is sufficient to prove that among the generic diffeomorphisms
studied in this section (i.e. the Moser generic diffeomorphisms, see
\ref{moserg_def}), the ones which don't have any periodic open annulus
without periodic points is generic. For such a map the interior of the
sets $E_q$, $q\geq1$, is empty as is the interior of the union $E$. We
deduce the second part of the proposition from the inclusion
$I\cup\partial E\subset \partial H$.

Thus it suffices to prove that there is a generic set of
diffeomorphisms in ${\rm Diff}_{\omega}^r(S^2)$ which have no periodic
open annulus which contains no periodic points.  Consider an open
connected and simply connected subset $V \subset S^2.$ We will first
show that for a dense open subset of ${\rm Diff}_{\omega}^r(S^2)$ the
set $V$ is not contained in a periodic open annulus which contains no
periodic points.

Consider the set of diffeomorphisms $\G_V \subset {\rm
Diff}_{\omega}^r(S^2)$ defined by the property that $F \in \G_V$ if
and only if there are two disjoint simple closed curves $C_1$ and
$C_2$ contained in one component of $\bigcup_{n \in Z} F^n(V)$ with
the property that there are hyperbolic periodic points of $F$ in each
of the three components of $S^2 \setminus (C_1 \cup C_2).$ It is clear
that if $F \in \G_V$ then $V$ is not contained in a periodic open
annulus which contains no periodic points since $C_1$ and $C_2$ would
have to be contained in such an annulus and this is impossible.  We will
show that $\G_V$ is dense and open.

The fact that $\G_V$ is open is clear since for any diffeomorphism $G$
close to $F$ we know $C_1 \cup C_2$ is contained in one component of
$\bigcup_{n \in Z} G^n(V)$ because $C_1 \cup C_2$ is compact.  Also
if $G$ is close to $F$ it will still have hyperbolic periodic points
in the three components of $S^2 \setminus (C_1 \cup C_2).$ So $\G_V$
is open.

To see that $\G_V$ is dense consider any $F_0 \in {\rm Diff}_{\omega}^r(S^2)$.
We can approximate it arbitrarily closely by a Moser generic diffeomorphism
$F$.  We will show either $F \in \G_V$ or it can be approximated by an
element of $\G_V.$  Consider $U = \bigcup_{n \in Z} F^n(V)$.  If
$U$ contains a periodic point of $F$ then it contains infinitely many
hyperbolic periodic points since for Moser generic diffeomorphisms
hyperbolic periodic points are not isolated and elliptic periodic points
are in the closure of hyperbolic periodic points.  It is thus clear that
if $U$ contains any periodic points we may choose simple closed curves
$C_1$ and $C_2$ contained in one component of $U$ such that each component
of $S^2 \setminus (C_1 \cup C_2)$ contains a hyperbolic periodic point.
Thus in this case $F \in \G_V.$

So we may assume $U$ contains no periodic points.  We choose a
periodic component $U_0$ of $U$, say of period $q$, so $F^q(U_0) =
U_0.$ If the complement of $U_0$ has three or more components
containing periodic points we can, as before, choose $C_1$ and $C_2$
showing $F \in \G_V.$ It is not possible for the complement of $U_0$
to have only one component containing periodic points, since the complement
of that component would be a simply connected $F^q$ invariant domain 
with no periodic points.  

All that remains is the case that the complement of $U_0$ has exactly
two components containing periodic points.  We choose a simple closed
curve $C$ in $U_0$ which separates these two components.  The
complement of these two components is an open $F^q$ invariant annulus
$W$ containing $U_0$ and in which $C$ is essential.

We choose a small annular neighborhood $U_1$ of $C$ in $U_0$.
Perturbing $F$ by the time $t$ map of an area preserving $C^\infty$
flow $\phi_t$ supported in the interior of $U_1$ we may continuously
vary the mean rotation number of $F^q \circ \phi_t = (F \circ
\phi_t)^q$ on the annulus $W$. However a result of \cite{F3} says that
if this number is rational there must exist a periodic point.  Clearly
if $t>0$ is small the perturbation is $C^r$ close to $F$.  One point on
the new periodic orbit must lie in the interior of $U_1$ since that is
where the perturbation is supported.  By a slight further perturbation
we obtain $G$ close to $F$ with a hyperbolic periodic point in the interior of 
$U_1$.  Choosing $C_1$ and $C_2$ to be the boundary components of 
$U_1$ we observe that $G \in \G_V$.

Thus we have shown that $\G_V$ is open and dense.  If we now choose
a countable basis $\{V_n\}$ for the topology of $S^2$ consisting of open
disks and let 
\[
\G = \bigcap_{n \ge 0} \G_{V_n}
\]
is the desired dense $G_{\delta}$.
\end{proof}

\end{document}